\documentstyle[12pt,amssymb]{article}
\def\Label{\label}

\long\def\ig#1{\relax}
\ig{Thanks to Roberto Minio for this def'n.  Compare the def'n of
\comment in AMSTeX.}

\newcount \coefa
\newcount \coefb
\newcount \coefc
\newdimen\tempdimen
\newdimen\xlen
\newdimen\ylen

\newcount\tempcounta
\newcount\tempcountb
\newcount\tempcountc
\newcount\tempcountd
\newcount\tempcounte
\newcount\tempcountf
\newcount\xext
\newcount\yext
\newcount\xoff
\newcount\yoff
\newcount\gap%
\newcount\arrowtypea
\newcount\arrowtypeb
\newcount\arrowtypec
\newcount\arrowtyped
\newcount\arrowtypee
\newcount\height
\newcount\width
\newcount\xpos
\newcount\ypos
\newcount\run
\newcount\rise
\newcount\arrowlength
\newcount\halflength
\newcount\arrowtype
\newsavebox{\tempboxa}%
\newsavebox{\tempboxb}%
\newsavebox{\tempboxc}%

\catcode`@=11 
\def\settoheight#1#2{\setbox\@tempboxa\hbox{#2}#1\ht\@tempboxa\relax}%
\def\settodepth#1#2{\setbox\@tempboxa\hbox{#2}#1\dp\@tempboxa\relax}%
\let\ifnextchar=\@ifnextchar
\catcode`@=12 

\def\putbox(#1,#2)#3{\put(#1,#2){\makebox(0,0){#3}}}

\def\setsqparms[#1`#2`#3`#4;#5`#6]{%
\settripairparms[#1`#2`#3`#4`1;#6]%
\width #5
}

\def\settriparms[#1`#2`#3;#4]{\settripairparms[#1`#2`#3`1`1;#4]}%

\def\settripairparms[#1`#2`#3`#4`#5;#6]{%
\arrowtypea #1
\arrowtypeb #2
\arrowtypec #3
\arrowtyped #4
\arrowtypee #5
\height #6
\width #6
}

\def\mvector(#1,#2)#3{
\put(0,0){\vector(#1,#2){#3}}%
\put(0,0){\vector(#1,#2){30}}%
}
\def\evector(#1,#2)#3{{
\arrowlength #3
\put(0,0){\vector(#1,#2){\arrowlength}}%
\advance \arrowlength by-30
\put(0,0){\vector(#1,#2){\arrowlength}}%
}}

\def\horsize#1#2{%
\settowidth{\tempdimen}{$#2$}%
#1=\tempdimen
\divide #1 by\unitlength
}

\def\vertsize#1#2{%
\settoheight{\tempdimen}{$#2$}%
#1=\tempdimen
\settodepth{\tempdimen}{$#2$}%
\advance #1 by\tempdimen
\divide #1 by\unitlength
}

\def\vertadjust[#1`#2`#3]{%
\vertsize{\tempcounta}{#1}%
\vertsize{\tempcountb}{#2}%
\ifnum \tempcounta<\tempcountb \tempcounta=\tempcountb \fi
\divide\tempcounta by2
\vertsize{\tempcountb}{#3}%
\ifnum \tempcountb>0 \advance \tempcountb by20 \fi
\ifnum \tempcounta<\tempcountb \tempcounta=\tempcountb \fi
}

\def\horadjust[#1`#2`#3]{%
\horsize{\tempcounta}{#1}%
\horsize{\tempcountb}{#2}%
\ifnum \tempcounta<\tempcountb \tempcounta=\tempcountb \fi
\divide\tempcounta by20
\horsize{\tempcountb}{#3}%
\ifnum \tempcountb>0 \advance \tempcountb by60 \fi
\ifnum \tempcounta<\tempcountb \tempcounta=\tempcountb \fi
}

\ig{ In this procedure, #1 is the paramater that sticks out all the way,
#2 sticks out the least and #3 is a label sticking out half way.  #4 is
the amount of the offset.}

\def\sladjust[#1`#2`#3]#4{%
\tempcountc=#4
\horsize{\tempcounta}{#1}%
\divide \tempcounta by2
\horsize{\tempcountb}{#2}%
\divide \tempcountb by2
\advance \tempcountb by-\tempcountc
\ifnum \tempcounta<\tempcountb \tempcounta=\tempcountb\fi
\divide \tempcountc by2
\horsize{\tempcountb}{#3}%
\advance \tempcountb by-\tempcountc
\ifnum \tempcountb>0 \advance \tempcountb by80\fi
\ifnum \tempcounta<\tempcountb \tempcounta=\tempcountb\fi
\advance\tempcounta by20
}

\def\putvector(#1,#2)(#3,#4)#5#6{{%
\xpos=#1
\ypos=#2
\run=#3
\rise=#4
\arrowlength=#5
\arrowtype=#6
\ifnum \arrowtype<0
    \ifnum \run=0
        \advance \ypos by-\arrowlength
    \else
        \tempcounta \arrowlength
        \multiply \tempcounta by\rise
        \divide \tempcounta by\run
        \ifnum\run>0
            \advance \xpos by\arrowlength
            \advance \ypos by\tempcounta
        \else
            \advance \xpos by-\arrowlength
            \advance \ypos by-\tempcounta
        \fi
    \fi
    \multiply \arrowtype by-1
    \multiply \rise by-1
    \multiply \run by-1
\fi
\ifnum \arrowtype=1
    \put(\xpos,\ypos){\vector(\run,\rise){\arrowlength}}%
\else\ifnum \arrowtype=2
    \put(\xpos,\ypos){\mvector(\run,\rise)\arrowlength}%
\else\ifnum\arrowtype=3
    \put(\xpos,\ypos){\evector(\run,\rise){\arrowlength}}%
\fi\fi\fi
}}

\def\bfig{\begin{picture}(\xext,\yext)(\xoff,\yoff)}
\def\efig{\end{picture}}

\def\putsplitvector(#1,#2)#3#4{
\xpos #1
\ypos #2
\arrowtype #4
\halflength #3
\arrowlength #3
\gap 140
\advance \halflength by-\gap
\divide \halflength by2
\ifnum \arrowtype=1
    \put(\xpos,\ypos){\line(0,-1){\halflength}}%
    \advance\ypos by-\halflength
    \advance\ypos by-\gap
    \put(\xpos,\ypos){\vector(0,-1){\halflength}}%
\else\ifnum \arrowtype=2
    \put(\xpos,\ypos){\line(0,-1)\halflength}%
    \put(\xpos,\ypos){\vector(0,-1)3}%
    \advance\ypos by-\halflength
    \advance\ypos by-\gap
    \put(\xpos,\ypos){\vector(0,-1){\halflength}}%
\else\ifnum\arrowtype=3
    \put(\xpos,\ypos){\line(0,-1)\halflength}%
    \advance\ypos by-\halflength
    \advance\ypos by-\gap
    \put(\xpos,\ypos){\evector(0,-1){\halflength}}%
\else\ifnum \arrowtype=-1
    \advance \ypos by-\arrowlength
    \put(\xpos,\ypos){\line(0,1){\halflength}}%
    \advance\ypos by\halflength
    \advance\ypos by\gap
    \put(\xpos,\ypos){\vector(0,1){\halflength}}%
\else\ifnum \arrowtype=-2
    \advance \ypos by-\arrowlength
    \put(\xpos,\ypos){\line(0,1)\halflength}%
    \put(\xpos,\ypos){\vector(0,1)3}%
    \advance\ypos by\halflength
    \advance\ypos by\gap
    \put(\xpos,\ypos){\vector(0,1){\halflength}}%
\else\ifnum\arrowtype=-3
    \advance \ypos by-\arrowlength
    \put(\xpos,\ypos){\line(0,1)\halflength}%
    \advance\ypos by\halflength
    \advance\ypos by\gap
    \put(\xpos,\ypos){\evector(0,1){\halflength}}%
\fi\fi\fi\fi\fi\fi
}

\def\setpos(#1,#2){\xpos=#1 \ypos#2}

\def\putmorphism(#1)(#2,#3)[#4`#5`#6]#7#8#9{{%
\run #2
\rise #3
\ifnum\rise=0
  \puthmorphism(#1)[#4`#5`#6]{#7}{#8}{#9}%
\else\ifnum\run=0
  \putvmorphism(#1)[#4`#5`#6]{#7}{#8}{#9}%
\else
\setpos(#1)%
\arrowlength #7
\arrowtype #8
\ifnum\run=0
\else\ifnum\rise=0
\else
\ifnum\run>0
    \coefa=1
\else
   \coefa=-1
\fi
\ifnum\arrowtype>0
   \coefb=0
   \coefc=-1
\else
   \coefb=\coefa
   \coefc=1
   \arrowtype=-\arrowtype
\fi
\width=2
\multiply \width by\run
\divide \width by\rise
\ifnum \width<0  \width=-\width\fi
\advance\width by60
\if l#9 \width=-\width\fi
\putbox(\xpos,\ypos){$#4$}
{\multiply \coefa by\arrowlength
\advance\xpos by\coefa
\multiply \coefa by\rise
\divide \coefa by\run
\advance \ypos by\coefa
\putbox(\xpos,\ypos){$#5$} }%
{\multiply \coefa by\arrowlength
\divide \coefa by2
\advance \xpos by\coefa
\advance \xpos by\width
\multiply \coefa by\rise
\divide \coefa by\run
\advance \ypos by\coefa
\if l#9%
   \put(\xpos,\ypos){\makebox(0,0)[r]{$#6$}}%
\else\if r#9%
   \put(\xpos,\ypos){\makebox(0,0)[l]{$#6$}}%
\fi\fi }%
{\multiply \rise by-\coefc
\multiply \run by-\coefc
\multiply \coefb by\arrowlength
\advance \xpos by\coefb
\multiply \coefb by\rise
\divide \coefb by\run
\advance \ypos by\coefb
\multiply \coefc by70
\advance \ypos by\coefc
\multiply \coefc by\run
\divide \coefc by\rise
\advance \xpos by\coefc
\multiply \coefa by140
\multiply \coefa by\run
\divide \coefa by\rise
\advance \arrowlength by\coefa
\ifnum \arrowtype=1
   \put(\xpos,\ypos){\vector(\run,\rise){\arrowlength}}%
\else\ifnum\arrowtype=2
   \put(\xpos,\ypos){\mvector(\run,\rise){\arrowlength}}%
\else\ifnum\arrowtype=3
   \put(\xpos,\ypos){\evector(\run,\rise){\arrowlength}}%
\fi\fi\fi}%
\fi\fi
\fi\fi}}

\def\puthmorphism(#1,#2)[#3`#4`#5]#6#7#8{{%
\xpos #1
\ypos #2
\width #6
\arrowlength #6
\putbox(\xpos,\ypos){$#3$\vphantom{$#4$}}%
{\advance \xpos by\arrowlength
\putbox(\xpos,\ypos){\vphantom{$#3$}$#4$}}%
\horsize{\tempcounta}{#3}%
\horsize{\tempcountb}{#4}%
\divide \tempcounta by2
\divide \tempcountb by2
\advance \tempcounta by30
\advance \tempcountb by30
\advance \xpos by\tempcounta
\advance \arrowlength by-\tempcounta
\advance \arrowlength by-\tempcountb
\putvector(\xpos,\ypos)(1,0){\arrowlength}{#7}%
\divide \arrowlength by2
\advance \xpos by\arrowlength
\vertsize{\tempcounta}{#5}%
\divide\tempcounta by2
\advance \tempcounta by20
\if a#8 %
   \advance \ypos by\tempcounta
   \put(\xpos,\ypos){\makebox(0,0){$#5$}}%
\else
   \advance \ypos by-\tempcounta
   \put(\xpos,\ypos){\makebox(0,0){$#5$}}%
\fi
}}

\def\putvmorphism(#1,#2)[#3`#4`#5]#6#7#8{{%
\xpos #1
\ypos #2
\arrowlength #6
\arrowtype #7
\settowidth{\xlen}{$#5$}%
\putbox(\xpos,\ypos){$#3$}%
{\advance \ypos by-\arrowlength
\putbox(\xpos,\ypos){$#4$}}%
{\advance\arrowlength by-140
\advance \ypos by-70
\ifdim\xlen>0pt
   \if m#8%
      \putsplitvector(\xpos,\ypos){\arrowlength}{\arrowtype}%
   \else
      \putvector(\xpos,\ypos)(0,-1){\arrowlength}{\arrowtype}%
   \fi
\else
   \putvector(\xpos,\ypos)(0,-1){\arrowlength}{\arrowtype}%
\fi}%
\ifdim\xlen>0pt
   \divide \arrowlength by2
   \advance\ypos by-\arrowlength
   \if l#8%
      \advance \xpos by-40
      \put(\xpos,\ypos){\makebox(0,0)[r]{$#5$}}%
   \else\if r#8%
      \advance \xpos by40
      \put(\xpos,\ypos){\makebox(0,0)[l]{$#5$}}%
   \else
      \putbox(\xpos,\ypos){$#5$}%
   \fi\fi
\fi
}}

\def\topadjust[#1`#2`#3]{%
\yoff=10
\vertadjust[#1`#2`{#3}]%
\advance \yext by\tempcounta
\advance \yext by 10
}

\def\botadjust[#1`#2`#3]{%
\vertadjust[#1`#2`{#3}]%
\advance \yext by\tempcounta
\advance \yoff by-\tempcounta
}

\def\leftadjust[#1`#2`#3]{%
\xoff=0
\horadjust[#1`#2`{#3}]%
\advance \xext by\tempcounta
\advance \xoff by-\tempcounta
}

\def\rightadjust[#1`#2`#3]{%
\horadjust[#1`#2`{#3}]%
\advance \xext by\tempcounta
}

\def\rightsladjust[#1`#2`#3]{%
\sladjust[#1`#2`{#3}]{\width}%
\advance \xext by\tempcounta
}

\def\leftsladjust[#1`#2`#3]{%
\xoff=0
\sladjust[#1`#2`{#3}]{\width}%
\advance \xext by\tempcounta
\advance \xoff by-\tempcounta
}

\def\adjust[#1`#2;#3`#4;#5`#6;#7`#8]{%
\topadjust[#1``{#2}]
\leftadjust[#3``{#4}]
\rightadjust[#5``{#6}]
\botadjust[#7``{#8}]}

\def\putsquare(#1)[#2`#3`#4`#5;#6`#7`#8`#9]{%
\setpos(#1)
\puthmorphism(\xpos,\ypos)[#4`#5`{#9}]{\width}{\arrowtyped}b%
\advance\ypos by \height
\puthmorphism(\xpos,\ypos)[#2`#3`{#6}]{\width}{\arrowtypea}a%
\putvmorphism(\xpos,\ypos)[``{#7}]{\height}{\arrowtypeb}l%
\advance\xpos by \width
\putvmorphism(\xpos,\ypos)[``{#8}]{\height}{\arrowtypec}r%
}

\def\square[#1`#2`#3`#4;#5`#6`#7`#8]{{
\xext=\width                              
\yext=\height                             
\topadjust[#1`#2`{#5}]
\botadjust[#3`#4`{#8}]
\leftadjust[#1`#3`{#6}]
\rightadjust[#2`#4`{#7}]
\begin{picture}(\xext,\yext)(\xoff,\yoff)
\putsquare(0,0)[#1`#2`#3`#4;#5`#6`#7`{#8}]
\end{picture}
}}

\def\putptriangle(#1,#2)[#3`#4`#5;#6`#7`#8]{%
\xpos=#1 \ypos=#2
\advance\ypos by \height
\puthmorphism(\xpos,\ypos)[#3`#4`{#6}]{\height}{\arrowtypea}a%
\putvmorphism(\xpos,\ypos)[`#5`{#7}]{\height}{\arrowtypeb}l%
\advance\xpos by\height
\putmorphism(\xpos,\ypos)(-1,-1)[``{#8}]{\height}{\arrowtypec}r%
}

\def\ptriangle[#1`#2`#3;#4`#5`#6]{{
\width=\height                         
\xext=\width                           
\yext=\width                           
\topadjust[#1`#2`{#4}]
\botadjust[#3``]
\leftadjust[#1`#3`{#5}]
\rightsladjust[#2`#3`{#6}]
\begin{picture}(\xext,\yext)(\xoff,\yoff)
\putptriangle(0,0)[#1`#2`#3;#4`#5`{#6}]%
\end{picture}%
}}

\def\putqtriangle(#1,#2)[#3`#4`#5;#6`#7`#8]{%
\xpos=#1 \ypos=#2
\advance\ypos by\height
\puthmorphism(\xpos,\ypos)[#3`#4`{#6}]{\height}{\arrowtypea}a%
\putmorphism(\xpos,\ypos)(1,-1)[``{#7}]{\height}{\arrowtypeb}l%
\advance\xpos by\height
\putvmorphism(\xpos,\ypos)[`#5`{#8}]{\height}{\arrowtypec}r%
}

\def\qtriangle[#1`#2`#3;#4`#5`#6]{{
\width=\height                         
\xext=\width                           
\yext=\height                          
\topadjust[#1`#2`{#4}]
\botadjust[#3``]
\leftsladjust[#1`#3`{#5}]
\rightadjust[#2`#3`{#6}]
\begin{picture}(\xext,\yext)(\xoff,\yoff)
\putqtriangle(0,0)[#1`#2`#3;#4`#5`{#6}]%
\end{picture}%
}}

\def\putdtriangle(#1,#2)[#3`#4`#5;#6`#7`#8]{%
\xpos=#1 \ypos=#2
\puthmorphism(\xpos,\ypos)[#4`#5`{#8}]{\height}{\arrowtypec}b%
\advance\xpos by \height \advance\ypos by\height
\putmorphism(\xpos,\ypos)(-1,-1)[``{#6}]{\height}{\arrowtypea}l%
\putvmorphism(\xpos,\ypos)[#3``{#7}]{\height}{\arrowtypeb}r%
}

\def\dtriangle[#1`#2`#3;#4`#5`#6]{{
\width=\height                         
\xext=\width                           
\yext=\height                          
\topadjust[#1``]
\botadjust[#2`#3`{#6}]
\leftsladjust[#2`#1`{#4}]
\rightadjust[#1`#3`{#5}]
\begin{picture}(\xext,\yext)(\xoff,\yoff)
\putdtriangle(0,0)[#1`#2`#3;#4`#5`{#6}]%
\end{picture}%
}}

\def\putbtriangle(#1,#2)[#3`#4`#5;#6`#7`#8]{%
\xpos=#1 \ypos=#2
\puthmorphism(\xpos,\ypos)[#4`#5`{#8}]{\height}{\arrowtypec}b%
\advance\ypos by\height
\putmorphism(\xpos,\ypos)(1,-1)[``{#7}]{\height}{\arrowtypeb}r%
\putvmorphism(\xpos,\ypos)[#3``{#6}]{\height}{\arrowtypea}l%
}

\def\btriangle[#1`#2`#3;#4`#5`#6]{{
\width=\height                         
\xext=\width                           
\yext=\height                          
\topadjust[#1``]
\botadjust[#2`#3`{#6}]
\leftadjust[#1`#2`{#4}]
\rightsladjust[#3`#1`{#5}]
\begin{picture}(\xext,\yext)(\xoff,\yoff)
\putbtriangle(0,0)[#1`#2`#3;#4`#5`{#6}]%
\end{picture}%
}}

\def\putAtriangle(#1,#2)[#3`#4`#5;#6`#7`#8]{%
\xpos=#1 \ypos=#2
{\multiply \height by2
\puthmorphism(\xpos,\ypos)[#4`#5`{#8}]{\height}{\arrowtypec}b}%
\advance\xpos by\height \advance\ypos by\height
\putmorphism(\xpos,\ypos)(-1,-1)[#3``{#6}]{\height}{\arrowtypea}l%
\putmorphism(\xpos,\ypos)(1,-1)[``{#7}]{\height}{\arrowtypeb}r%
}

\def\Atriangle[#1`#2`#3;#4`#5`#6]{{
\width=\height                         
\xext=\width                           
\yext=\height                          
\topadjust[#1``]
\botadjust[#2`#3`{#6}]
\multiply \xext by2 
\leftsladjust[#2`#1`{#4}]
\rightsladjust[#3`#1`{#5}]
\begin{picture}(\xext,\yext)(\xoff,\yoff)%
\putAtriangle(0,0)[#1`#2`#3;#4`#5`{#6}]%
\end{picture}%
}}

\def\putAtrianglepair(#1,#2)[#3]{\xpos=#1 \ypos=#2%
\putAtrianglepairx[#3]}
\def\putAtrianglepairx[#1`#2`#3`#4;#5`#6`#7`#8`#9]{%
\puthmorphism(\xpos,\ypos)[#2`#3`{#8}]{\height}{\arrowtyped}b%
\advance\xpos by\height
\puthmorphism(\xpos,\ypos)[\phantom{#3}`#4`{#9}]{\height}{\arrowtypee}b%
\advance\ypos by\height
\putmorphism(\xpos,\ypos)(-1,-1)[#1``{#5}]{\height}{\arrowtypea}l%
\putvmorphism(\xpos,\ypos)[``{#6}]{\height}{\arrowtypeb}m%
\putmorphism(\xpos,\ypos)(1,-1)[``{#7}]{\height}{\arrowtypec}r%
}

\def\Atrianglepair[#1`#2`#3`#4;#5`#6`#7`#8`#9]{{%
\width=\height
\xext=\width
\yext=\height
\topadjust[#1``]%
\vertadjust[#2`#3`{#8}]
\tempcountd=\tempcounta                       
\vertadjust[#3`#4`{#9}]
\ifnum\tempcounta<\tempcountd                 
\tempcounta=\tempcountd\fi                    
\advance \yext by\tempcounta                  
\advance \yoff by-\tempcounta                 
\multiply \xext by2 
\leftsladjust[#2`#1`{#5}]
\rightsladjust[#4`#1`{#7}]%
\begin{picture}(\xext,\yext)(\xoff,\yoff)%
\putAtrianglepair(0,0)[#1`#2`#3`#4;#5`#6`#7`#8`{#9}]%
\end{picture}%
}}

\def\putVtriangle(#1,#2)[#3`#4`#5;#6`#7`#8]{%
\xpos=#1 \ypos=#2
\advance\ypos by\height
{\multiply\height by2
\puthmorphism(\xpos,\ypos)[#3`#4`{#6}]{\height}{\arrowtypea}a}%
\putmorphism(\xpos,\ypos)(1,-1)[`#5`{#7}]{\height}{\arrowtypeb}l%
\advance\xpos by\height
\advance\xpos by\height
\putmorphism(\xpos,\ypos)(-1,-1)[``{#8}]{\height}{\arrowtypec}r%
}

\def\Vtriangle[#1`#2`#3;#4`#5`#6]{{
\width=\height                         
\xext=\width                           
\yext=\height                          
\topadjust[#1`#2`{#4}]
\botadjust[#3``]
\multiply \xext by2 
\leftsladjust[#1`#3`{#5}]
\rightsladjust[#2`#3`{#6}]
\begin{picture}(\xext,\yext)(\xoff,\yoff)%
\putVtriangle(0,0)[#1`#2`#3;#4`#5`{#6}]%
\end{picture}%
}}

\def\putVtrianglepair(#1,#2)[#3]{\xpos=#1 \ypos=#2%
\putVtrianglepairx[#3]}
\def\putVtrianglepairx[#1`#2`#3`#4;#5`#6`#7`#8`#9]{%
\advance\ypos by\height
\putmorphism(\xpos,\ypos)(1,-1)[`#4`{#7}]{\height}{\arrowtypec}l%
\puthmorphism(\xpos,\ypos)[#1`#2`{#5}]{\height}{\arrowtypea}a%
\advance\xpos by\height
\puthmorphism(\xpos,\ypos)[\phantom{#2}`#3`{#6}]{\height}{\arrowtypeb}a%
\putvmorphism(\xpos,\ypos)[``{#8}]{\height}{\arrowtyped}m%
\advance\xpos by\height
\putmorphism(\xpos,\ypos)(-1,-1)[``{#9}]{\height}{\arrowtypee}r%
}

\def\Vtrianglepair[#1`#2`#3`#4;#5`#6`#7`#8`#9]{{%
\xoff=0
\yoff=2 
\xext=\height                  
\width=\height                 
\yext=\height                  
\vertadjust[#1`#2`{#5}]
\tempcountd=\tempcounta        
\vertadjust[#2`#3`{#6}]
\ifnum\tempcounta<\tempcountd
\tempcounta=\tempcountd\fi
\advance \yext by\tempcounta
\botadjust[#4``]%
\multiply \xext by2
\leftsladjust[#1`#4`{#7}]%
\rightsladjust[#3`#4`{#9}]%
\begin{picture}(\xext,\yext)(\xoff,\yoff)%
\putVtrianglepair(0,0)[#1`#2`#3`#4;#5`#6`#7`#8`{#9}]%
\end{picture}%
}}

\def\putCtriangle(#1,#2)[#3`#4`#5;#6`#7`#8]{%
\xpos=#1 \ypos=#2
\advance\ypos by\height
\putmorphism(\xpos,\ypos)(1,-1)[``{#8}]{\height}{\arrowtypec}l%
\advance\xpos by\height
\advance\ypos by\height
\putmorphism(\xpos,\ypos)(-1,-1)[#3`#4`{#6}]{\height}{\arrowtypea}l%
{\multiply\height by 2
\putvmorphism(\xpos,\ypos)[`#5`{#7}]{\height}{\arrowtypeb}r}%
}

\def\Ctriangle[#1`#2`#3;#4`#5`#6]{{
\width=\height                          
\xext=\width                            
\yext=\height                           
\multiply \yext by2 
\topadjust[#1``]
\botadjust[#3``]
\sladjust[#2`#1`{#4}]{\width}
\tempcountd=\tempcounta                 
\sladjust[#2`#3`{#6}]{\width}
\ifnum \tempcounta<\tempcountd          
\tempcounta=\tempcountd\fi              
\advance \xext by\tempcounta            
\advance \xoff by-\tempcounta           
\rightadjust[#1`#3`{#5}]
\begin{picture}(\xext,\yext)(\xoff,\yoff)%
\putCtriangle(0,0)[#1`#2`#3;#4`#5`{#6}]%
\end{picture}%
}}

\def\putDtriangle(#1,#2)[#3`#4`#5;#6`#7`#8]{%
\xpos=#1 \ypos=#2
\advance\xpos by\height \advance\ypos by\height
\putmorphism(\xpos,\ypos)(-1,-1)[``{#8}]{\height}{\arrowtypec}r%
\advance\xpos by-\height \advance\ypos by\height
\putmorphism(\xpos,\ypos)(1,-1)[`#4`{#7}]{\height}{\arrowtypeb}r%
{\multiply\height by 2
\putvmorphism(\xpos,\ypos)[#3`#5`{#6}]{\height}{\arrowtypea}l}%
}

\def\Dtriangle[#1`#2`#3;#4`#5`#6]{{
\width=\height                         
\xext=\height                          
\yext=\height                          
\multiply \yext by2 
\topadjust[#1``]
\botadjust[#3``]
\leftadjust[#1`#3`{#4}]
\sladjust[#2`#1`{#4}]{\height}
\tempcountd=\tempcountd                
\sladjust[#2`#3`{#6}]{\height}
\ifnum \tempcounta<\tempcountd         
\tempcounta=\tempcountd\fi             
\advance \xext by\tempcounta           
\begin{picture}(\xext,\yext)(\xoff,\yoff)
\putDtriangle(0,0)[#1`#2`#3;#4`#5`{#6}]%
\end{picture}%
}}

\def\setrecparms[#1`#2]{\width=#1 \height=#2}%
\def\recurse[#1`#2`#3`#4;#5`#6`#7`#8`#9]{{%
\settowidth{\tempdimen}{#1}
\ifdim\tempdimen=0pt
  \savebox{\tempboxa}{\hbox{#2}}%
  \savebox{\tempboxb}{\hbox{#4}}%
  \savebox{\tempboxc}{\hbox{#7}}%
\else
  \savebox{\tempboxa}{\hbox{$\hbox{#1}\times\hbox{#2}$}}%
  \savebox{\tempboxb}{\hbox{$\hbox{#1}\times\hbox{#4}$}}%
  \savebox{\tempboxc}{\hbox{$\hbox{#1}\times\hbox{#7}$}}%
\fi
\tempcounte=\height
\divide\tempcounte by 2
\tempcountf=\tempcounte
\advance\tempcountf by \width
\xext=\tempcountf \yext=\height
\topadjust[#2`\usebox{\tempboxa}`{#5}]%
\botadjust[#4`\usebox{\tempboxb}`{#9}]%
\sladjust[#3`#2`{#6}]{\tempcounte}%
\tempcountd=\tempcounta
\sladjust[#3`#4`{#8}]{\tempcounte}%
\ifnum \tempcounta<\tempcountd
\tempcounta=\tempcountd\fi
\advance \xext by\tempcounta
\advance \xoff by-\tempcounta
\rightadjust[\usebox{\tempboxa}`\usebox{\tempboxb}`\usebox{\tempboxc}]%
\bfig
{\settriparms[-1`1`1;\tempcounte]%
\putCtriangle(0,0)[`#3`;#6`#7`{#8}]}%
\arrowtypea=-1 \arrowtypeb=0 \arrowtypec=1 \arrowtyped=-1
\putsquare(\tempcounte,0)[#2`\usebox{\tempboxa}`#4`\usebox{\tempboxb};%
#5``\usebox{\tempboxc}`#9]%
\efig
}}

\newtheorem{definition}{Definition $\!\!$}[section]
\newtheorem{prop}[definition]{Proposition $\!\!$}
\newtheorem{lem}[definition]{Lemma $\!\!$}
\newtheorem{corollary}[definition]{Corollary $\!\!$}
\newtheorem{theorem}[definition]{Theorem $\!\!$}
\newtheorem{example}[definition]{Example $\!\!$}
\newtheorem{remark}[definition]{Remark $\!\!$}
\newtheorem{con}[definition]{Conjecture $\!\!$}

\newcommand{\nc}[2]{\newcommand{#1}{#2}}
\newcommand{\rnc}[2]{\renewcommand{#1}{#2}}
\nc{\bpr}{\begin{prop}}
\nc{\bth}{\begin{theorem}}
\nc{\ble}{\begin{lem}}
\nc{\bco}{\begin{corollary}}
\nc{\bre}{\begin{remark}}
\nc{\bex}{\begin{example}}
\nc{\bde}{\begin{definition}}
\nc{\bcon}{\begin{con}}
\nc{\econ}{\end{con}}
\nc{\ede}{\end{definition}}
\nc{\eco}{\end{corollary}}
\nc{\ere}{\hfill\mbox{$\Diamond$}\end{remark}}
\nc{\eex}{\end{example}}
\nc{\esa}{\end{satz}}
\nc{\epr}{\end{prop}}
\nc{\ethe}{\end{theorem}}
\nc{\ele}{\end{lem}}
\nc{\epf}{\hfill\mbox{$\Box$}\\~\\}
\nc{\beq}{\begin{equation}}
\nc{\eeq}{\end{equation}}
\nc{\ot}{\otimes}
\nc{\lra}{\longrightarrow}
\nc{\ci}{\circ}
\nc{\ba}{\begin{array}}
\nc{\ea}{\end{array}}
\nc{\bea}{\begin{eqnarray}}
\nc{\eea}{\end{eqnarray}}
\rnc{\[}{\beq}
\rnc{\]}{\eeq}
\nc{\bpf}{{\it Proof:}~ }
\def\Z{{\Bbb Z}}
\def\N{{\Bbb N}}
\def\R{{\Bbb R}}
\def\C{{\Bbb C}}
\def\Sh{{\frak{S}}}
\rnc{\phi}{\varphi}

\nc{\qpb}{quantum principal bundle}
\nc{\te}{\!\ot\!}
\nc{\pf}{\mbox{$P\!\sb F$}}
\nc{\pn}{\mbox{$P\!\sb\nu$}}
\nc{\bmlp}{\mbox{\boldmath$\left(\right.$}}
\nc{\bmrp}{\mbox{\boldmath$\left.\right)$}}
\rnc{\phi}{\mbox{$\varphi$}}
\nc{\LAblp}{\mbox{\LARGE\boldmath$($}}
\nc{\LAbrp}{\mbox{\LARGE\boldmath$)$}}
\nc{\Lblp}{\mbox{\Large\boldmath$($}}
\nc{\Lbrp}{\mbox{\Large\boldmath$)$}}
\nc{\lblp}{\mbox{\large\boldmath$($}}
\nc{\lbrp}{\mbox{\large\boldmath$)$}}
\nc{\blp}{\mbox{\boldmath$($}}
\nc{\brp}{\mbox{\boldmath$)$}}
\nc{\LAlp}{\mbox{\LARGE $($}}
\nc{\LArp}{\mbox{\LARGE $)$}}
\nc{\Llp}{\mbox{\Large $($}}
\nc{\Lrp}{\mbox{\Large $)$}}
\nc{\llp}{\mbox{\large $($}}
\nc{\lrp}{\mbox{\large $)$}}
\nc{\lbc}{\mbox{\Large\boldmath$,$}}
\nc{\lc}{\mbox{\Large$,$}}
\nc{\Lall}{\mbox{\Large$\forall\;$}}
\nc{\bc}{\mbox{\boldmath$,$}}
\rnc{\epsilon}{\varepsilon}
\rnc{\ker}{\mbox{\rm Ker}}
\nc{\ra}{\rightarrow}
\nc{\la}{\triangleright}

\nc{\cc}{\!\ci\!}
\nc{\T}{\mbox{\sf T}}
\nc{\can}{\mbox{\em\sf T}\!\sb R}
\nc{\cnl}{$\mbox{\sf T}\!\sb R$}

\nc{\M}{\mbox{Map}}
\rnc{\to}{\mapsto}
\nc{\imp}{\Rightarrow}
\rnc{\iff}{\Leftrightarrow}
\nc{\ob}{\mbox{$\Omega\sp{1}\! (\! B)$}}
\nc{\op}{\mbox{$\Omega\sp{1}\! (\! P)$}}
\nc{\oa}{\mbox{$\Omega\sp{1}\! (\! A)$}}
\nc{\inc}{\mbox{$\,\subseteq\;$}}

\nc{\spp}{\mbox{${\cal S}{\cal P}(P)$}}
\nc{\dr}{\mbox{$\Delta_{R}$}}
\nc{\dsr}{\mbox{$\Delta_{\cal R}$}}
\nc{\m}{\mbox{m}}
\nc{\0}{\sb{(0)}}
\nc{\1}{\sb{(1)}}
\nc{\2}{\sb{(2)}}
\nc{\3}{\sb{(3)}}
\nc{\4}{\sb{(4)}}
\nc{\5}{\sb{(5)}}
\nc{\6}{\sb{(6)}}
\nc{\7}{\sb{(7)}}
\nc{\hsp}{\hspace*}
\nc{\al}{\mbox{$u$}}
\nc{\bet}{\mbox{$\beta$}}
\nc{\ha}{\mbox{$\alpha$}}
\nc{\hb}{\mbox{$\beta$}}
\nc{\hg}{\mbox{$\gamma$}}
\nc{\hd}{\mbox{$\delta$}}
\nc{\he}{\mbox{$\varepsilon$}}
\nc{\hz}{\mbox{$\zeta$}}
\nc{\hs}{\mbox{$\sigma$}}
\nc{\hk}{\mbox{$\kappa$}}
\nc{\hm}{\mbox{$\mu$}}
\nc{\hn}{\mbox{$\nu$}}
\nc{\hl}{\mbox{$\lambda$}}
\nc{\hG}{\mbox{$\Gamma$}}
\nc{\hD}{\mbox{$\Delta$}}

\nc{\Th}{\mbox{$\Theta$}}
\nc{\ho}{\mbox{$\omega$}}
\nc{\hO}{\mbox{$\Omega$}}
\nc{\hp}{\mbox{$\pi$}}
\nc{\hP}{\mbox{$\Pi$}}

\nc{\as}{\mbox{$A(S^3\sb s)$}}
\nc{\bs}{\mbox{$A(S^2\sb s)$}}
\nc{\slq}{\mbox{$A(SL\sb q(2))$}}
\nc{\fr}{\mbox{$Fr\llp A(SL(2,\IC))\lrp$}}
\nc{\slc}{\mbox{$A(SL(2,\IC))$}}
\nc{\af}{\mbox{$A(F)$}}
\rnc{\widetilde}{\tilde}
\nc{\suq}{\mbox{$A(SU_q(2))$}}
\nc{\asq}{\mbox{$A(S_q^2)$}}
\nc{\tasq}{\mbox{$\widetilde{A}(S_q^2)$}}

\def\st{\stackrel}

\newcommand{\fa}{\forall}

\def\IC{{\Bbb C}}
\def\IZ{{\Bbb Z}}
\def\IN{{\Bbb N}}
\def\IR{{\Bbb R}}

\def\ra{\rightarrow}

\def\<{\langle}
\def\>{\rangle}

\begin{document}
\baselineskip16pt
\parskip=.5\baselineskip

\title{\vspace*{-1.5cm}\LARGE\bf
QUANTUM REAL PROJECTIVE SPACE, DISC AND SPHERE\\
{\normalsize\it Dedicated to the memory of Stanis\l aw Zakrzewski.}\\
}
\author{
\vspace*{-1mm}\large\sc 
Piotr M.~Hajac\\
\vspace*{-2mm}\normalsize 
Mathematical Institute, Polish Academy of Sciences\\
\vspace*{-1mm}\normalsize 
ul.~\'Sniadeckich 8, Warsaw, 00--950 Poland\\
\vspace*{-1mm}\normalsize 
and\\
\vspace*{-2mm}\normalsize 
Department of Mathematical Methods in Physics\\
\vspace*{-1mm}\normalsize 
Warsaw University, ul. Ho\.{z}a 74, Warsaw, 00-682 Poland\\
\large
http://www.fuw.edu.pl/$\!\!\!\!\!\!\widetilde{\phantom{mmm}}\!\!\!\!\!\!$pmh
\normalsize 
\and
\vspace*{-1mm}\large\sc
Rainer Matthes\\
\vspace*{-2mm}\normalsize
Max Planck Institute for Mathematics in the Sciences\\
\vspace*{-1mm}\normalsize 
Inselstr. 22--26, D-04103 Leipzig, Germany\\
\vspace*{-1mm}\normalsize 
and\\
\vspace*{-2mm}\normalsize 
Institute of Theoretical Physics, Leipzig University\\
\vspace*{-1mm}\normalsize 
Augustusplatz 10/11, D-04109 Leipzig, Germany\\
\large
e-mail: matthes@itp.uni-leipzig.de
\and
\vspace*{-1mm}\large\sc 
Wojciech Szyma\'nski\\
\vspace*{-2mm}\normalsize 
School of Mathematical and Physical Sciences, University of Newcastle\\
\vspace*{-1mm}\normalsize 
Callaghan, NSW 2308, Australia\\
\large 
e-mail: wojciech@frey.newcastle.edu.au
}
\date{\normalsize }
\maketitle

{\bf Abstract.}
We define the $C^*$-algebra of quantum real projective space $\R P_q^2$,
classify its irreducible representations and compute its $K$-theory.
We also show that the $q$-disc of Klimek-Lesniewski  can be 
obtained as a non-Galois $\Z_2$-quotient of the 
equator Podle\'s  quantum sphere. 
On the way, we provide the Cartesian coordinates for 
all Podle\'s quantum spheres and
determine an explicit form of isomorphisms between the
 $C^*$-algebras of the equilateral spheres and the $C^*$-algebra 
of the equator one.


\section{Introduction}

Classical spheres can be constructed by gluing two discs along their
boundaries. Since an open disc is homeomorphic to $\R^2$, this fact is 
reflected
in the following short exact sequence of $C^*$-algebras
of continuous functions (vanishing at infinity where appropriate):
\beq
\Label{cs}
0\lra C_0(\R^2)\oplus C_0(\R^2)\lra C(S^2)\lra C(S^1)\lra 0.
\eeq
On the other hand, 
one can obtain a disc $D^2$ as the quotient of a sphere under
the $\Z_2$-action given by the reflection with respect to the equator plane.
 Two copies of an open disc
collapse to one copy, and we have the short exact sequence
\beq
\Label{cd}
0\lra C_0(\R^2)\lra C(D^2)\lra C(S^1)\lra 0.
\eeq
 Similarly, real
projective space $\IR P^2$ can be constructed from the antipodal action of
$\Z_2$ on the two-sphere. As for $D^2$, removing $S^1$ from $\R P^2$
also leaves an open disc, and again we have the short exact sequence
\beq
\Label{crp}
0\lra C_0(\R^2)\lra C(\R P^2)\lra C(S^1)\lra 0.
\eeq
 The aim of this paper is to present the noncommutative
geometry of a $q$-deformation of the aforementioned geometric setting.
(This deformation is unique under some assumptions.)
It turns out that the $q$-deformation changes $C_0(\R^2)$ in the above
short exact sequences into the ideal ${\cal K}$ of compact operators
(see (\ref{qs}), (\ref{kdi}), (\ref{qrp})). Therefore, since $C_0(\R^2)$
and $\cal K$ behave in a similar way in $K$-theory, it is not surprising
that the $K$-groups of these $q$-deformed surfaces coincide with the
respective $K$-groups of their classical counterparts.
Since $D^2$ has a boundary and $\R P^2$ is non-orientable, we hope that
the study of their $q$-analogues will help one to understand the concept
of a boundary and orientability in the general noncommutative setting.

Deformations of $SL(2,\C)$ were studied in depth and classified
\cite{dl90,w-sl91,wz94}. The choice of the compact $*$-structure and
the requirement of the existence of the $C^*$-norm lead then to the celebrated
deformation of $SU(2)$, which we denote by $SU_q(2)$.
(The literature on this quantum group motivating and 
treating it from many different points of view is vast. E.g., see \cite{ks97}
for references.) 
Subsequently, the study of quantum
homogeneous spaces of $SU_q(2)$ leads to the classification of
quantum spheres~\cite{p-p87}. (See \cite{s-a91} for the Poisson aspects.)
On the other hand,
motivated by the Poisson geometry, noncommutative deformations of the unit disc
were constructed in \cite{kl92,kl93}. Gluings of quantum discs
which produce quantum spheres were studied in~\cite{s-a91,mnw91,cm00}. 
Finally, quantum real projective space $\IR P^2_q$ was defined in \cite{h-pm96}
within the framework of the Hopf-Galois theory to exemplify the concept of
strong connections on quantum principal bundles (cf.~\cite[Example 2.13]{dgh}). 
It was obtained as the quantum quotient space from the antipodal $\Z_2$-action
on the Podle\'s equator sphere. This action was already discovered 
 in~\cite{p-p87}, and is
the only possible  $\Z_2$-action on quantum spheres compatible with the actions
of~$SU_q(2)$
(see above Section~6 therein).

In this paper, we continue along these lines. We begin in Section~2 by
reviewing the relevant known results on quantum spheres ($C^*$-representations,
$K$-theory). Then we provide the Cartesian coordinates
 and compute an explicit form of the 
$C^*$-isomorphisms between the $C^*$-algebra of the equator quantum sphere
and the $C^*$-algebras of the equilateral Podle\'s spheres ($c\in(0,\infty)$).
We also show that these isomorphisms commute with the $U(1)$-actions inherited
from the actions of $SU_q(2)$ on quantum spheres. In Section~3, we prove that
the $q$-disc of Klimek-Lesniewski can be obtained as a noncommutative quotient
of the equator quantum sphere by an appropriate $\Z_2$-action. More precisely,
first we show that the polynomial algebra of the $q$-disc is a fixed-point
subalgebra of the polynomial algebra of the equator quantum sphere under a
non-Galois $\Z_2$-action. Then we extend this construction to the equilateral
quantum spheres by employing the aforementioned $C^*$-isomorphisms. Since
these isomorphisms are non-polynomial, we handle the equilateral spheres only
on the $C^*$-level. We complete this section by recalling the topological
$K$-theory of the $q$-disc. The paper ends with Section~4 where we define
the $C^*$-algebra of quantum $\IR P^2$, study its representations, and compute
the $K$-theory. Similarly to the quantum disc case, this $C^*$-algebra is
obtained as a $\Z_2$-action fixed-point subalgebra of the $C^*$-algebra of 
the equator quantum sphere.
For both the quantum disc and $\IR P^2_q$ cases, 
we show  that the $\Z_2$-actions are compatible with the above-mentioned
actions of $U(1)$.

Throughout the paper we use the jargon of Noncommutative Geometry referring
to quantum spaces as objects dual  to noncommutative algebras
in the sense of the Gelfand-Naimark correspondence
between spaces and function algebras. 
The unadorned tensor product means the completed (spatial) tensor product
when placed between $C^*$-algebras, and the algebraic tensor product over
$\C$ otherwise. The algebras are assumed to be associative and over $\C$.
They are also unital unless the contrary is obvious from the context.
By $P(\mbox{quantum space})$ we denote the polynomial algebra of a
quantum space, and by $C(\mbox{quantum space})$ the corresponding $C^*$-algebra. 
In this paper, the $C^*$-completion ($C^*$-closure) of a $*$-algebra always
means the completion with respect to the supremum norm over all
$*$-representations in bounded operators.

\section{Quantum spheres}

\bde[\cite{p-p87}]
The $C^*$-algebra $C(S^2_{q\infty})$ of the quantum sphere $S^2_{q\infty}$,
$q\in\IR$, $0<|q|<1$, is defined as the $C^*$-closure of the $*$-algebra
$P(S^2_{q\infty}):={\IC}\<A,B\>/I_{q\infty}$, where $I_{q\infty}$ 
is the (two-sided) $*$-ideal
in the free $*$-algebra ${\IC}\<A,B\>$ generated by the relations
\beq\Label{ab}
A^*=A, ~~~ BA=q^2AB, 
\eeq
\beq\Label{b*b}
B^*B=-A^2+I, ~~~
BB^*=-q^4A^2+I .
\eeq
The $C^*$-algebra $C(S^2_{qc})$ of the quantum sphere $S^2_{qc}$, $c\in 
[0,\infty)$ is defined analogously, with  (\ref{b*b}) 
replaced by
\beq\Label{b*bc}
B_c^*B_c=A_c-A_c^2+cI,~~~
B_cB_c^*=q^2A_c-q^4A_c^2+cI .
\eeq
\ede
The irreducible $*$-representations of the quantum spheres are determined in
\cite{p-p87}. Let us denote by $\pi^c_\pm$ and by $\pi^c_\theta$ the infinite
dimensional and one-dimensional representations of $C(S^2_{qc})$ ($c\in[0,\infty]$),
respectively. (In the $c=\infty$ case, in agreement with the notation for generators
in the above definition, we write $\pi_\pm$ and $\pi_\theta$ instead of
 $\pi^\infty_\pm$ and $\pi^\infty_\theta$, respectively.)
The complete list of the irreducible $*$-representations of $C(S^2_{q\infty})$  
is given by 
\beq
\pi_\theta(A)=0,~~~ \pi_\theta(B)=e^{i\theta},~~~ \theta\in [0,2\pi),
\eeq
and 
\beq\Label{rs2inf}
\pi_\pm(A)e_k=\pm q^{2k}e_k, ~~~
\pi_\pm(B)e_k=(1-q^{4k})^{1/2}e_{k-1},~~~
\pi_{\pm}(B)e_0=0.
\eeq
Here $\{e_k\}_{k\in\N}$ is an orthonormal basis of a Hilbert space.
Similarly, the irreducible $*$-representations of 
$C(S^2_{qc})$, $c\in(0,\infty)$, are defined by
\beq\Label{theta}
\pi_{\theta}^c(A_c)=0,~~~ \pi_{\theta}^c(B_c)=c^{1/2}e^{i\theta},
~~~ \theta\in [0,2\pi),
\eeq
and 
\beq\Label{rs2ca}
\pi_{\pm}^c(A_c)e_k=\lambda_\pm q^{2k}e_k,\;\;\;
\pi_{\pm}^c(B_c)e_k={c_{\pm}(k)}^{1/2}e_{k-1},\;\;\;
\pi_{\pm}^c(B_c)e_0=0,
\eeq
\beq\Label{lamb}
\mbox{where }~\lambda_{\pm}=\frac{1}{2}\pm(c+\frac{1}{4})^{1/2},\;\;\;
c_{\pm}(k)=\lambda_{\pm}q^{2k}-(\lambda_\pm q^{2k})^2+c.
\eeq
The direct sums $\pi_+^c\oplus\pi_-^c$, $0<c\leq\infty$, are
faithful representations. 
The representations $\pi^c_\pm$ can be considered as embeddings 
of quantum discs onto
the northern and  southern hemisphere, respectively, whereas
the one-dimensional representations are the classical points (forming a circle).
For $c=\infty$, the classical points are
symmetric with respect to the hemispheres, i.e., they form the equator. 
With $c$ decreasing, the circle of classical points shrinks to a pole. 
Thus, in the limit case $c=0$, we can think of a quantum sphere as a quantum disc
whose (classical) boundary is glued to a point. 
For $c=0$, the formulas (\ref{theta})-(\ref{lamb}) still define $*$-representations.
Now, however, $\pi^0_\theta$ coincide for all $\theta$, and  $\pi^0_\theta$ and
$\pi^0_+$ are the only irreducible representations. The representation $\pi^0_-$ 
becomes trivial, and $\pi^0_+$ becomes faithful.
The cases $c=0$, $0<c<\infty$ and $c=\infty$ are referred to as the standard,
equilateral and equator quantum sphere, respectively.

To make the aforementioned geometric picture explicit, we need to find
the {\em Cartesian coordinates for quantum spheres}. More precisely,
we need to define self-adjoint generators $x,y,z$ of 
$P(S^2_{qc}),~c\in[0,\infty],$
which satisfy
$x^2+y^2+z^2=1$.
Note first that dividing (\ref{b*bc}) by $c$ and rescaling the generators
by $c^{-1/2}$ would lead to the formulas whose limit with $c\ra\infty$
 would be (\ref{b*b}). To include also the $c=0$ case, let us rescale the
generators by $(1+\sqrt{c})^{-1}$, i.e., 
\beq\Label{scale}
\widetilde{A}_c:=\frac{A_c}{1+\sqrt{c}},~~~
\widetilde{B}_c:=\frac{B_c}{1+\sqrt{c}}.
\eeq
Now, from (\ref{b*bc}), we have
\beq\Label{b*bt}
\widetilde{B}_c^*\widetilde{B}_c=\frac{\widetilde{A}_c}{1+\sqrt{c}}
-\widetilde{A}_c^2+\frac{c}{(1+\sqrt{c})^2}I,~~~
\widetilde{B}_c\widetilde{B}_c^*=\frac{q^2\widetilde{A}_c}{1+\sqrt{c}}
-q^4\widetilde{A}_c^2+\frac{c}{(1+\sqrt{c})^2}I.
\eeq
The relations (\ref{ab}) remain unchanged, that is,
$\widetilde{A}_c^*=\widetilde{A}_c,~\widetilde{B}_c\widetilde{A}_c=
q^2\widetilde{A}_c\widetilde{B}_c$. 
Contrary to $A_c$ and $B_c$, the generators $\widetilde{A}_c$ and
$\widetilde{B}_c$ have limits with $c\ra\infty$ when thought of as
elements of $P(SU_q(2))$. Indeed, remembering the definition of
$A_c,B_c$ \cite[pp.196,200]{p-p87} 
in terms of the spin 1 representation 
\beq
D_1:=\left(\ba{ccc}
\hd^2 & -(1+q^2)\hd\hg & -q\hg^2\\
-q^{-1}\hb\hg & 1-(q+q^{-1})\hb\hg & \ha\hg\\
-q^{-1}\hb^2 & -(q+q^{-1})\hb\ha & \ha^2
\ea\right)
\eeq
of $SU_q(2)$ 
(with $\ha,\hb,\hg,\hd$ being the generators of the algebra $P(SU_q(2))$),
 we can write
\beq\Label{matrix}
\left(\widetilde{B}_c^*,\widetilde{A}_c,\widetilde{B}_c\right)
=\left(\mbox{$\frac{\sqrt{c}}{1+\sqrt{c}}$},\mbox{$\frac{1}{1+\sqrt{c}}$},
\mbox{$\frac{\sqrt{c}}{1+\sqrt{c}}$}\right)
D_1
\left(\ba{ccc}
1&0&0\\
0&-(1+q^2)^{-1}&0\\
0&0&1
\ea\right)
+\left(0,\mbox{$\frac{1}{(1+\sqrt{c})(1+q^2)}$},0\right).
\eeq
It is clear now that the tilded generators are well-defined also
for $c=\infty$. Since the relations among the tilded generators
become for $c=\infty$ the relations among $A$ and $B$, we can write
$\widetilde{A}_\infty:=A$, $\widetilde{B}_\infty:=B$.
Thus, we have a uniform description
of quantum spheres for all $c\in[0,\infty]$. 
(See \cite[Section 6]{bm00} for a uniform
parameterization of Podle\'s spheres by the unit interval $[0,1]$.)
Remembering the geometrical meaning of $\widetilde{A}_c,
\widetilde{B}_c$ (see \cite[pp.196,200,201]{p-p87}), we put
\beq\Label{cart1}
x=iQ_x(\widetilde{B}_c-\widetilde{B}_c^*),
~~~
y=Q_y(\widetilde{B}_c+\widetilde{B}_c^*),
~~~
z={Q}_z(\widetilde{A}_c-\widetilde{a}_0),~~~
c\in[0,\infty ].
\eeq
Here $Q_x$, $Q_y$, $Q_z$ and $\widetilde{a}_0$
are real-valued functions of $q$ and~$c$, so that $x,y,z$ are evidently
self-adjoint.
The condition $x^2+y^2+z^2=1$  and the linear independence
of the monomials 
$A_c^kB_c^l,~A_c^m{B_c^*}^n,\; 
k,l,m,n \in {\IN}, n>0$ \cite[p.116]{p-p89} imply that 
$Q_x^2=Q_y^2$. Let us put 
$Q_h:=|Q_x|=|Q_y|$. Then the sphere equation boils down to
\beq\Label{boil}
2Q_h^2(\widetilde{B}_c\widetilde{B}_c^*
+\widetilde{B}_c^*\widetilde{B}_c)
+Q_z^2(\widetilde{A}_c-\widetilde{a}_0)^2=1.
\eeq
Plugging in (\ref{b*bt}) to the above formula yields
\beq
2Q_h^2\left(-(1+q^4)\widetilde{A}_c^2+\frac{1+q^2}{1+\sqrt{c}}\widetilde{A}_c+
\frac{2c}{(1+\sqrt{c})^2}\right)
+Q_z^2(\widetilde{A}_c-\widetilde{a}_0)^2=1.
\eeq
Employing again the linear independence of the monomials $A_c^k$, 
one can compute:
\beq
Q_h=-\mbox{$\frac{\sqrt{2(1+q^4)}}{1+q^2}$}(1+\sqrt{c})z_\infty,~~~
|Q_z|=-\mbox{$\frac{2(1+q^4)}{1+q^2}$}(1+\sqrt{c})z_\infty,~~~
\widetilde{a}_0=\mbox{$\frac{1+q^2}{2(1+q^4)}$}(1+\sqrt{c})^{-1},
\eeq
where $z_\infty:=-\left(8c\frac{1+q^4}{(1+q^2)^2}+1\right)^{-1/2}$.
(The meaning of this number will shortly become clear.)
Let us choose $Q_x=Q_h=Q_y$, $Q_z=|Q_z|$. The formulas (\ref{cart1})
read now:
\bea\Label{cart2}
&&
x=-i\mbox{$\frac{\sqrt{2(1+q^4)}}{1+q^2}$}(1+\sqrt{c})z_\infty
(\widetilde{B}_c-\widetilde{B}_c^*),
\nonumber\\ &&
y=-\mbox{$\frac{\sqrt{2(1+q^4)}}{1+q^2}$}(1+\sqrt{c})z_\infty
(\widetilde{B}_c+\widetilde{B}_c^*),
\nonumber\\ &&
z=-\mbox{$\frac{2(1+q^4)}{1+q^2}$}(1+\sqrt{c})z_\infty
\widetilde{A}_c+z_\infty.
\eea
(Observe that $z_\infty|_{c=\infty}=0$ and 
$((1+\sqrt{c})z_\infty)|_{c=\infty}=-\frac{1+q^2}{2\sqrt{2(1+q^4)}}$.)
The eigenvalues ${z}^\pm_k$ of $\pi^c_\pm({z})$ 
 are given by
\beq
\pi^c_\pm({z})e_k=
\left(z_\infty-z_\infty\hl_\pm\frac{2(1+q^4)}{1+q^2}q^{2k}\right) e_k,
~~~ c\in[0,\infty].
\eeq
(Note that $(z_\infty\hl_\pm)|_{c=\infty}=\mp\frac{1+q^2}{2\sqrt{2(1+q^4)}}$.)
It is evident that $\lim_{k\ra\infty}z^\pm_k=z_\infty$. Since we also have
$
\pi_\theta^c({z})={z}_\infty,
$
we can say that the eigenvalues of $\pi^c_\pm({z})$ converge (from both sides) to
the circle of classical points (space of one-dimensional representations) given
by $\pi_\theta^c$. For $c=\infty$ we have $z_\infty=0$, so that the circle is the
equator, whereas for $c= 0$ the circle shrinks to the south pole ($z_\infty=-1$).
Finally, let us remark that, as $(1+\sqrt{c})z_\infty\neq 0$ for
any $c\in[0,\infty]$, the equations (\ref{cart2}) can be solved for 
$\widetilde{A}_c,\widetilde{B}_c,\widetilde{B}_c^*$, and consequently $x,y,z$
generate the algebra $P(S^2_{qc})$. Since they are also self-adjoint and satisfy
$x^2+y^2+z^2=1$, we call them the Cartesian coordinates of quantum spheres.  

We recall from \cite{s-a91} that $\pi_+^c\oplus\pi_-^c$
is for all $c\in(0,\infty]$ a $C^*$-isomorphism of $C(S^2_{qc})$ onto
$C^*({\Sh})\oplus_\sigma C^*({\Sh})$.
Here $C^*({\Sh})$
is the $C^*$-algebra of the one-sided shift (Toeplitz algebra).
It is the $C^*$-algebra generated by
the shift operator ${\Sh}e_i=e_{i+1}$, where $\{e_i\}_{i\in\IN}$ is 
an orthonormal
basis of a Hilbert space. The map $\sigma:C^*({\Sh})\ra C(S^1)$ is 
the so-called symbol map defined by ${\Sh}\mapsto u$, where
$u$ is the unitary generator of
$C(S^1)$. The algebra
$C^*({\Sh})\oplus_\sigma 
C^*({\Sh})$ is
defined as the gluing of two copies of $C^*(\Sh)$ via $\sigma$, i.e.,
\[
C^*({\Sh})\oplus_\sigma 
C^*({\Sh}):=\{(a_1,a_2)\in C^*({\Sh})\oplus 
C^*({\Sh})|\sigma(a_1)=\sigma(a_2)\}.
\]
Let
\beq
\chi_c:=(\pi_+\oplus\pi_-)^{-1}\ci(\pi^c_+\oplus\pi^c_-):
C(S^2_{qc})\lra C(S^2_{q\infty})
\eeq
 be
 the isomorphism composed from the isomorphisms
$\pi_+\oplus\pi_-:C(S^2_{q\infty})\ra C^*(\Sh)\oplus_\sigma C^*(\Sh)$ and
$\pi_+^c\oplus\pi_-^c:C(S^2_{qc})\ra C^*(\Sh)\oplus_\sigma C^*(\Sh)$.
An explicit form of the isomorphisms $\chi_c$ is given by:
\bpr\Label{FG}
Let $\eta_c(t):=\sqrt{t-t^2+c}$ (cf.~(\ref{lamb})), and let
$F_c$ and $G_c$ be functions given by
\[F_c(x):=\left\{\begin{array}{cl}\lambda_+x&\mbox{for}~~~~~0\leq x\leq 1\\
-\lambda_-x&\mbox{for}~~-1\leq
x< 0, \end{array}\right.
\]
\[
G_c(x):=(1-q^4x^2)^{-1/2}\left\{\begin{array}{cl}\eta_c(q^2\hl_+x)
&\mbox{for}~~~~~0\leq x\leq 1\\
\eta_c(-q^2\hl_-x)&\mbox{for}~~-1\leq x< 0. \end{array}\right.
\]
Then
$\chi_c(A_c)=F_c(A)$ and $\chi_c(B_c)=G_c(A)B$.
\epr
\bpf
First, note that, since $\pi_\pm(A)$ is diagonal and $F_c$ and $G_c$ are continuous
functions defined on the spectrum of $\pi_\pm(A)$, the operators $F_c(\pi_\pm(A))$
and $G_c(\pi_\pm(A))$ make sense and are easily computable.
Subsequently,
notice that $\chi_c(A_c)=F_c(A)$ if and only if $\pi^c_\pm(A_c)=F_c(\pi_\pm(A))$.
To verify the latter equality, we check that 
\[
\pi_{\pm}^c(A_c)e_k
=\lambda_{\pm}q^{2k}e_k
=F_c(\pi_\pm(A))e_k.
\]
Similarly, to verify $\chi_c(B_c)=G_c(A)B$,
we observe that 
$G_c(\pm q^{2k})=\frac{c_\pm(k+1)^{1/2}}{(1-q^{4(k+1)})^{1/2}}$ and compute
\[
\pi^c_\pm(B_c)e_k
=\pi_\pm(G_c(A)B)e_k
=G_c(\pi_\pm(A))\pi_\pm(B)e_k,
\]
which proves the proposition.
\epf
The above proposition shows that the isomorphisms $\chi_c$ are of non-polynomial
nature. Therefore we suspect that:
\bcon\Label{cc}
The polynomial $*$-algebras of the quantum spheres $S^2_{qc}$ and $S^2_{qc'}$
are non-isomorphic for $c\neq c'$.
\econ

Our next step is to consider the compatibility of the isomorphisms $\chi_c$
with the actions of $U(1)$ inherited from the actions of $SU_q(2)$ on quantum
spheres.
Let $c\in(0,\infty]$, and let $\delta:C(S^2_{qc})\ra C(S^2_{qc})\ot C(U(1))$ be
 the right coaction obtained from the coaction
$\Delta_R:C(S^2_{qc})\ra
C(S^2_{qc})\ot C(SU_q(2))$ \cite[p.194]{p-p87} with the help of the natural map
$C(SU_q(2))\ra C(U(1))$. Explicitly, we have
\beq\Label{delta}
\delta(A_c)=A_c\otimes 1, ~~\delta(B_c)=B_c\otimes u^2.
\eeq
Here $u$ is the unitary generator of $C(U(1))$. 
Since $\hd:C(S^2_{qc})\ra C(S^2_{qc})\ot C(U(1))$
is a $*$-homomorphism, it is continuous. Therefore, 
identifying $C(S^2_{qc})\ot C(U(1))$ with $C(U(1),C(S^2_{qc}))$
(continuous functions on $U(1)$ with values in $C(S^2_{qc})$; see 
\cite[Proposition T.5.21]{w-ne93}), we obtain, for any $g\in U(1)$, a continuous map
\beq\Label{delta2}
\hd_g:C(S^2_{qc})\lra C(S^2_{qc}),~~~ \hd_g(a):=\hd(a)(g).
\eeq
This defines an action of $U(1)$ on $S^2_{qc}$. It follows from the continuity of 
$\hd$ and the fact that
the image of $\hd$ is in the continuous functions from $U(1)$ to
$C(S^2_{qc})$ that each $\hd_g$ is a $C^*$-automorphism
of $C(S^2_{qc})$.
Contrary to the action of $SU_q(2)$,
the action of $U(1)$ on $S^2_{qc}$ is compatible with the quantum 
``homeomorphisms'' among the spheres $S^2_{qc}$, i.e.,
$(\chi^{-1}_{c'}\ci\chi_c)\ci\hd_g=\delta_g\ci(\chi^{-1}_{c'}\ci
\chi_c).$ This follows from:

\bpr\Label{com}
$\forall g\in U(1):~~\delta_g\ci\chi_c=\chi_c\ci\hd_g.$
\epr
\bpf
Since both $\chi_c$ and $\hd_g$ are $C^*$-isomorphisms, 
it suffices to check this equality
on generators. It follows from (\ref{delta}) that $\hd_g(A_c)=A_c$
 and $\hd_g(B_c)=g^2B_c$. Taking advantage of Proposition~\ref{FG},
one can compute:
\[
(\delta_g\ci\chi_c)(A_c)
=\delta_g(F_c(A))
=F_c(\delta_g(A))
=F_c(A)=\chi_c(A_c)
=(\chi_c\ci\hd_g)(A_c),
\]
\[
(\delta_g\ci\chi_c)(B_c)
=\delta_g(G_c(A)B)
=G_c(\delta_g(A))\delta_g(B)
=g^2G_c(A)B=g^2\chi_c(B_c)
=(\chi_c\ci\hd_g)(B_c).
\]
This proves the proposition.
\epf

For the sake of completeness (cf.~(\ref{cs}), (\ref{kdi}), (\ref{qrp})),
let us end  this 
section by recalling the topological K-theory of the quantum spheres. 
First, there is an exact sequence~\cite[Proposition 1.2]{s-a91}: 
\beq 
\Label{qs}
0\lra{\cal K}\oplus{\cal K}\lra C(S^2_{qc})
   \lra C(S^1)\lra 0, 
\eeq
where ${\cal K}$ is the ideal of compact operators. It induces the 6-term
exact sequence in K-theory, from which it follows 
that $K_0 (C(S^2_{qc}))\cong\IZ\oplus\IZ$, $K_1(C(S^2_{qc}))\cong 0$
\cite[Proposition 4.1]{mnw91}.

\section{Quantum disc}
\setcounter{equation}{0}

\bde[\cite{kl93}]
The $C^*$-algebra $C(D_q)$, $0<q < 1$, of the quantum disc $D_q$ is 
 the $C^*$-closure (obtained from $*$-representations in bounded operators)
of the algebra $P(D_q):={\IC}\<x,x^*\>/J_q$, 
where $J_q$ is the two-sided ideal
in the free algebra ${\IC}\<x,x^*\>$ generated by the relation
\beq\Label{disc}
x^*x-qxx^*=1-q.
\eeq
\ede
The goal of this section is to determine the relationship between the thus
defined quantum discs and the equator and equilateral quantum spheres
(cf.\ \cite[p.278]{nn94} and references therein).
The objects $D_q$ form
 a one-parameter sub-family of the two-parameter family of quantum
discs described in \cite{kl93}. Explicitly, the latter family is given by
$x^*x-qxx^*=1-q + \mu(xx^*-1)(x^*x-1)$. It is known 
(\cite[Proposition VI.1]{kl93}, 
\cite[Proposition 15]{cm00}, \cite[Theorem IV.7]{kl92}, \cite[p.222]{s-a91}) that
for $0\leq\mu<1-q$ and $q=1$, $0<\mu <1$, 
the quantum-disc $C^*$-algebras are all isomorphic to the Toeplitz
algebra (the $C^*$-algebra generated by the one-sided shift 
${\frak S}e_i=e_{i+1}$). Furthermore, we know from \cite{kl93} that
every irreducible {\em bounded} $*$-representation of $P(D_q)$ 
is unitarily equivalent 
to a one-dimensional representation $\pi_{\theta}$ defined by
\beq\Label{th}
\pi_{\theta}(x)=e^{i\theta},~~~ \pi_{\theta}(x^*)=e^{-i\theta},
~~ 0\leq\theta<2\pi,
\eeq
or an infinite dimensional representation $\pi$ given on
an orthonormal basis $\{e_i\}_{i\in\N}$ by the formulas
\beq \Label{pix}
\pi(x)e_i=(1-q^{i+1})^{1/2}~e_{i+1},~~~i\geq 0,
\eeq
\beq\Label{pix*}
\pi(x^*)e_i=\left\{\begin{array}{lr}0,&i=0,\\(1-q^i)^{1/2}~e_{i-1},&
i\geq 1.
\end{array}\right.
\eeq
The infinite dimensional representation 
$\pi$ is faithful \cite[p.14]{kl93}. Also, one can directly verify that
$\pi$ is faithful on
the polynomial algebra $P(D_q)$, so that $P(D_q)\inc C(D_q)$.
Finally, let us mention that there are also 
unbounded representations of the relation (\ref{disc}). 
They are given, e.g., in \cite[Section 5.2.6]{ks97}.

Now we are going to show that the above-defined $q$-disc can be
obtained by collapsing the equator quantum sphere.
In the classical case, the $\Z_2$-action on $S^2$ 
collapsing it to a disc is not free,
as it leaves the equator invariant. This entails that the map
\beq\Label{psi}
\psi:S^2\times \Z_2\ni(x,g)\longmapsto (x,xg)\in S^2\times S^2
\eeq
is not injective, whence $S^2\ra D^2$ is not a principal fibration. (The considered
$\Z_2$-action is not Galois.)
The $\Z_2$-action on quantum spheres that we are looking for should identify 
``points'' of the same type and leave the equator invariant. Therefore, since
the standard sphere contains only one classical point (pole), we exclude it
from our considerations. 
Our first step is to define the desired $\Z_2$-action on the polynomial
algebra of the equator quantum sphere. As in the classical case, we define it as
the reflection with respect to the equator plane, i.e., via the $*$-algebra
automorphism $r_1$ of $P(S^2_{q\infty})$ sending $B$ to itself, and $A$ to $-A$.
(It is immediate from the commutation relations of the equator quantum sphere 
that $r_1$ is well defined.) 
Dualizing the $\Z_2$-action $r_1$ on $S^2_{q\infty}$ we get the coaction
$\hD_1:P(S^2_{q\infty})\ra P(S^2_{q\infty})\ot P(\Z_2)$ making 
$P(S^2_{q\infty})$ a right
$P(\Z_2)$-comodule algebra. (See \cite[Section 1.6]{m-s93} for generalities.)
Explicitly, denoting by $\la$ the action of $\Z_2$ on $P(S^2_{q\infty})$, we have
\[
\hD_1(p)= (1\la p)\ot 1^*+((-1)\la p)\ot(-1)^*=
\frac{1}{2}\llp p\ot(1+\ha)+r_1(p)\ot(1-\ha)\lrp.
\]
Here $\{1^*,(-1)^*\}$ denotes the basis dual to the basis $\{1,-1\}$ of the group ring
$\C[\Z_2]$, and $\alpha(\pm 1)=\pm 1$.
The main claim of this section is that the quantum disc is a non-Galois
quotient of the equator quantum sphere:
\bth\Label{main} The polynomial algebra of the equator quantum sphere
is a non-Galois $\IZ_2$-extension of the polynomial algebra of the
quantum disc via the above defined action $r_1$, i.e.,

1. $P(D^2_{q^4})\cong P(S^2_{q\infty}/\Z_2):= 
\{a\in P(S^2_{q\infty})~|~r_1(a)=a\}$ ($\Z_2$-extension).

2. The canonical map
$
P(S^2_{q\infty})\ot_{P(D^2_{q^4})}P(S^2_{q\infty})\ni p\ot_{P(D^2_{q^4})}p'\to
p\hD_1(p')\in P(S^2_{q\infty})\ot P(\Z_2)
$
\phantom{is {\em not}} is {\em not} bijective. (The extension is not Galois.)
\ethe
\bpf
1. We know from \cite[p.116]{p-p89} that the monomials 
\[
A^kB^l,~A^m{B^*}^n,~ 
k,l,m,n \in {\IN}, n>0,
\]
 form a linear basis of $P(S^2_{q\infty})$. 
Since $r_1(A)=-A$ and $r_1(B)=B$,
taking advantage of the above basis, one can see that $r_1(a)=a$ if and only if
$a$ is a linear combination of basis monomials that have $A$ in the even
power. It follows now from (\ref{b*b}) that any $r_1$-invariant $a$
 is a polynomial
in $B$ and $B^*$. Thus, since every polynomial in $B$ and $B^*$
 is $r_1$-invariant,
$P(S^2_{q\infty}/\Z_2)$ is the $*$-subalgebra generated by $B$.
 On the other hand,
one can conclude from (\ref{b*b}) that
\beq\Label{bdisc}
BB^*-q^4B^*B=(1-q^4)1.
\eeq
This equation, together with (\ref{disc}), allows us to define a 
$*$-epimorphism 
\[
\phi:P(D_{q^4})\lra P(S^2_{q\infty}/\Z_2),
~~\phi(x)= B^*.
\]
 To complete the proof, we need
to show that $\phi$ is injective. It is immediate from formula (\ref{rs2inf})
that $\pi_\pm\ci\phi=\pi$, where $\pi$ is defined by
(\ref{pix})-(\ref{pix*}). Hence the injectivity of $\pi$
implies the injectivity of~$\phi$. 

2. It suffices to show that the map 
\[\hk:P(S^2_{q\infty})\ot P(S^2_{q\infty})\ni p\ot p'\to
p\hD_1(p')\in P(S^2_{q\infty})\ot P(\Z_2)\]
is not surjective. (The considered $\Z_2$-action on $S^2_{q\infty}$ is not free.)
Note first that in the classical case to detect the lack of surjectivity of 
the pullback
map $\psi^*$ (see (\ref{psi})), we can use the function $1\ot\ha$. 
Indeed, for any point $x$ on the equator, we have
\[\psi^*(\mbox{anything})(x,1)=(\mbox{anything})(x,x)=\psi^*(\mbox{anything})(x,-1),\]
whereas
$(1\ot\ha)(x,1)=1\neq-1=(1\ot\ha)(x,-1)$.
It turns out that $1\ot\ha$ also does the job in the quantum case.
Suppose that $1\ot\ha$ is in the image of $\hk$. Then there exists a tensor
$\sum_ip_i\ot p'_i$ such that
\[
\frac{1}{2}\sum_ip_i\llp p'_i\ot (1+\ha)+r_1(p'_i)\ot(1-\ha)\lrp=1\ot\ha.
\]
Evaluating the right tensorands on both sides at $1$ and $-1$ yields
\[
\left\{\begin{array}{l}\sum_ip_ip'_i=1\\ \sum_ip_ir_1(p'_i)=-1.\end{array}\right.
\]
Applying $\pi_\theta$ (see (\ref{th})) to these equations gives
\[
\left\{\begin{array}{l}\sum_i\pi_\theta(p_i)\pi_\theta(p'_i)=1\\ 
\sum_i\pi_\theta(p_i)
\pi_\theta(r_1(p'_i))= \sum_i\pi_\theta(p_i)
\pi_\theta(p'_i)=-1,\end{array}\right.
\]
which is the desired contradiction.
\epf

In order to define a $\Z_2$-action on the closure $C(S^2_{q\infty})$ of
$P(S^2_{q\infty})$,
note that the flip map 
\[
\tau:C({\Sh})\oplus_\sigma C({\Sh})\ra
C({\Sh})\oplus_\sigma C({\Sh}),~~~ \tau(a,b)=(b,a),
\]
 satisfies
$\tau\ci(\pi_+\oplus\pi_-)=(\pi_+\oplus\pi_-)\ci r_1$. It is therefore natural
to define the completion of $r_1$ to a $C^*$-algebra map by
\footnote{We owe this idea to S.L.\ Woronowicz.}
\beq\Label{br1}
\bar{r}_1:=(\pi_+\oplus\pi_-)^{-1}\ci\tau\ci(\pi_+\oplus\pi_-).
\eeq
\bpr\Label{co}
The $C^*$-subalgebra 
$C(S^2_{q\infty}/\Z_2):=\{a\in C(S^2_{q\infty})|\bar{r}_1(a)=a\}$
of $\Z_2$-invariants in $C(S^2_{q\infty})$ coincides
with the $C^*$-completion of $P(S^2_{q\infty}/\Z_2)$ inside $C(S^2_{q\infty})$,
and is isomorphic to
$C(D_{q^4})$.
\epr
\bpf
First let us argue that the map $\phi$ defined in the proof of Proposition 
\ref{main}
extends to a $C^*$-isomorphism of $C(D_{q^4})$ with the closure of
$P(S^2_{q\infty}/\Z_2)$ inside $C(S^2_{q\infty})$.
Since $\phi:P(D^2_{q^4})\ra P(S^2_{q\infty}/\Z_2)$ is an isomorphism of
$*$-algebras, the $*$-representations of $P(D^2_{q^4})$ can be turned to
$*$-representations of $P(S^2_{q\infty}/\Z_2)$, and vice-versa. We are to show
that $\phi$ determines a one-to-one correspondence between the $*$-representations
used to define the norm on $C(D^2_{q^4})$ and $C(S^2_{q\infty})$ respectively.
As every $*$-representation of $P(S^2_{q\infty})$ is bounded, it yields via
$\phi$ a {\em bounded} $*$-representation of $P(D^2_{q^4})$. On the other 
hand, every bounded $*$-representation of $P(D^2_{q^4})$ gives a $*$-representation 
of $P(S^2_{q\infty}/\Z_2)$ which can be extended to  $C(S^2_{q\infty})$.
Indeed, let $P(D^2_{q^4})\stackrel{\rho}{\ra} B(H)$ be such a representation.
Then $\widetilde{\rho}:=\rho\ci\phi^{-1}$ is a bounded $*$-representation of
$P(S^2_{q\infty}/\Z_2)$, and since $B$ and $B^*$ satisfy the disc relation
(\ref{bdisc}), it follows from \cite{kl93} that $||\widetilde{\rho}(B^*B)||=1$. 
As a consequence, $I-\widetilde{\rho}(B^*B)\geq 0$, and one can define 
$\widetilde{\rho}(A)=\sqrt{I-\widetilde{\rho}(B^*B)}$.
This gives the desired extension. 

To complete the proof, note first that,
since $\bar{r}_1$ is continuous, the $C^*$-closure of $P(S^2_{q\infty}/\Z_2)$
is contained in $C(S^2_{q\infty}/\Z_2)$. Thus  it only
remains to show that every $\Z_2$-invariant in $C(S^2_{q\infty})$
is in the closure of $P(S^2_{q\infty}/\Z_2)$. Let $a \in C(S^2_{q\infty}/\Z_2)$.
Then, again by the continuity of $\bar{r}_1$ and density of $P(S^2_{q\infty})$
in $C(S^2_{q\infty})$, $a=\lim_{n\rightarrow\infty}a_n$ with $a_n\in 
P(S^2_{q\infty})$ and
\[
a= \frac{1}{2}(id+\bar{r}_1)(a)
=\frac{1}{2}(id+\bar{r}_1)(\lim_{n\rightarrow\infty}a_n)
=\lim_{n\rightarrow\infty}\frac{1}{2}(id+r_1)(a_n).
\]
As $\frac{1}{2}(id+r_1)(a_n)\in P(S^2_{q\infty}/\Z_2)$ for any $n$, $a$
is in the closure of $P(S^2_{q\infty}/\Z_2)$, as claimed.
\epf
\vspace*{-10mm}
\bre\em
Since all quantum disc algebras $C(D_q)$, $0<q<1$, are isomorphic 
as $C^*$-algebras to the $C^*$-algebra of the one-sided shift, in the above
proposition $q^4$ can be replaced by any element of the 
interval $(0,1)$.
\ere

We extend the $\Z_2$-action to the equilateral quantum spheres by the formula
\beq\Label{r1c}
\bar{r}_1^c:=\chi_c^{-1}\ci\bar{r}_1\ci\chi_c.
\eeq
It is now evident that we have
\bco
The subalgebra $C(S^2_{qc}/\Z_2):=\{a\in C(S^2_{qc})~|~\bar{r}^c_1(a)=a\}$
of $\Z_2$-invariants of $C(S^2_{qc})$, $c\in(0,\infty]$, is isomorphic to the
$C^*$-algebra $C(D_{q^4})$ of the quantum disc.
\eco
Furthermore, it is clear from (\ref{br1}) and (\ref{r1c}) that
\beq
\bar{r}^c_1=(\pi^c_+\oplus\pi^c_-)^{-1}\ci \tau\ci(\pi^c_+\oplus\pi^c_-).
\eeq
Explicitly, the above equality reads
\footnote{We are grateful to S.~L.~Woronowicz for putting us on the track
of the reasoning below.}
\[
\pi_+^c(a)
=\pi_-^c(\bar{r}^c_1(a)),~~\pi_-^c(a)
=\pi_+^c(\bar{r}^c_1(a)),~~a\in C(S^2_{qc}).
\]
For $a=A_c$, these equations are solved by the formulas
\beq\Label{ra}
\bar{r}^c_1(A_c)= f_c(A_c),~~~
f_c(x)=\left\{\begin{array}{cc}
\frac{\lambda_-}{\lambda_+}x&x\geq 0\\
\frac{\lambda_+}{\lambda_-}x&x\leq 0.\end{array}\right.
\eeq
Note that this piecewise linear function $f_c$ can be replaced by 
any continuous function 
having the same
values as $f_c$ at the points $\lambda_\pm q^{2k}$. For $c\in(0,\infty)$,
among these functions there is no 
polynomial.
This shows that $\bar{r}^c_1$ does {\em not} leave $P(S^2_{qc})$ invariant. 
Furthermore, considering 
the  image of $B_c^*$ under $\pi_+^c\oplus\pi_-^c$, one finds the polar
decomposition $B_c^*=V_c|B_c^*|$ with 
\beq\Label{pol}
\pi_{\pm}^c(V_c)e_k=e_{k+1},~~~
\pi_{\pm}^c(|B_c^*|)e_k=c_\pm(k+1)^{1/2}e_k,~~~k\in\IN.
\eeq
Hence  $(\pi_+^c\oplus\pi_-^c)(V_c)=(\Sh,\Sh)$,
and $\bar{r}_1^c(V_c)=V_c$.
Recall that the spectrum of $(\pi_+^c\oplus\pi_-^c)(A_c)$ 
and $(\pi_+^c\oplus\pi_-^c)(|B_c^*|)$
is $\{0\}\cup\{\lambda_\pm q^{2k}~|~k\in \IN\}$ and 
$\{\sqrt{c}\}\cup\{\sqrt{c_\pm(k+1)}~|~k\in \IN\}$ respectively. 
On the other hand, 
$\pi_{\pm}^c(\bar{r}_1^c(|B_c^*|))e_k=c_\mp(k+1)^{1/2}e_k,$ $k\in\IN$.
One can directly check that
\beq\Label{rb*}
\bar{r}_1^c(|B_c^*|)
=g_c(A_c),~~ 
g_c(t)
:=\left\{\begin{array}{cc}\eta_c(\frac{\lambda_-}{\lambda_+}q^2t)&0\leq t\leq
\lambda_+\\\eta_c(\frac{\lambda_+}{\lambda_-}q^2t)&\lambda_-\leq t<0,
\end{array}\right.~~~
\eta_c(t)=\sqrt{t-t^2+c}.
\eeq
Clearly, $g_c$ can be replaced by any continuous function having at the points 
$\lambda_\pm q^{2k}$ values $c_\mp(k+1)^{1/2}$, for any $k$.
Note that we used $A_c$ instead of $|B_c^*|$ to obtain $\bar{r}_1^c(|B_c^*|)$ 
as a continuous function of a generator because the assignment 
$\sqrt{c_+(k+1)}\mapsto \sqrt{c_-(k+1)}$ does not give a function, as 
$k_1\neq k_2$ implies $c_-(k_1)\neq c_-(k_2)$, whereas it might happen that
$c_+(k_1)=c_+(k_2)$ for $k_1\neq k_2$.
Indeed, let $k_1$, $k_2$ be any two different positive
natural numbers. Then the equation $\sqrt{c_+(k_1)}=\sqrt{c_+(k_2)}$ is equivalent
to the equation $q^{2k_1}+q^{2k_2}=\lambda_+^{-1}$. Since $\lambda_+^{-1}\in (0,1)$,
there exists $q\in (0,1)$ solving this equality.
\bre\em
The formulas for $f_c, F_c, g_c, G_c$ are consistent with one another 
by construction. Nevertheless, it is entertaining
to verify this consistency in a direct manner. 
Taking into account Proposition~\ref{FG} and
formula (\ref{ra}), we obtain a sequence of equivalent equalities
\begin{eqnarray*}
\bar{r}_1^c(A_c)&=&(\chi_c^{-1}\ci\bar{r}_1\ci\chi_c)(A_c)\\
f_c(A_c)&=&(\chi_c^{-1}\ci\bar{r}_1)(F_c(A))\\
f_c(A_c)&=&\chi_c^{-1}(F_c(-A))\\
\pi_c(f_c(A_c))&=&\pi(F_c(-A))\\
f_c(\pi_c(A_c))&=&F_c(-\pi(A))\\
f_c(\pi^c_\pm(A_c))e_k&=&F_c(\pi_\mp(A))e_k.
\end{eqnarray*}
Recalling (\ref{rs2inf}) and (\ref{rs2ca}) one can see that the last equality is true.
Similarly, taking advantage of the polar decomposition $B_c^*=V_c|B_c^*|$,
$\bar{r}_1^c(V_c)=V_c$ and (\ref{rb*}), we get
\begin{eqnarray*}
\bar{r}_1^c(B^*_c)&=&(\chi_c^{-1}\ci\bar{r}_1\ci\chi_c)(B^*_c)\\
\bar{r}_1^c(V_c|B^*_c|)&=&(\chi_c^{-1}\ci\bar{r}_1)(\chi_c(B^*_c))\\
\bar{r}_1^c(V_c)\bar{r}_1^c(|B^*_c|)&=&(\chi_c^{-1}\ci\bar{r}_1)(B^*G_c(A))\\
V_cg_c(A_c)&=&\chi_c^{-1}(B^*G_c(-A))\\
\pi_c(V_cg_c(A_c))&=&\pi(B^*G_c(-A))\\
\pi_\pm^c(V_c)g_c(\pi_\pm^c(A_c))e_k&=&\pi_\pm(B^*)G_c(\pi_\mp(A))e_k.
\end{eqnarray*}
Remembering formulas (\ref{rs2inf}), (\ref{rs2ca}), (\ref{pol}), the last equality
is evident.
\ere

Next, 
let us consider the rotational invariance with respect to the South-North Pole axis
of the above-studied $\Z_2$-actions on quantum spheres.
The $U(1)$-action on $S^2_{q\infty}$ 
(see (\ref{delta2})) is given on generators by $\hd_g(A)=A$, $\hd_g(B)=g^2B$.
Therefore, one can infer from Proposition~\ref{com} that the $U(1)$-action on 
$S^2_{qc}$ ($c\in(0,\infty]$)
 and the reflection $\bar{r}_1^c$ are compatible:
\beq\Label{coma}
\hd_g\ci\bar{r}_1^c=\bar{r}_1^c\ci\hd_g.
\eeq
\bre\em\Label{comp}
It follows already from (\ref{coma}) that
$
\delta\ci \bar{r}_1^c=(\bar{r}_1^c\ot id)\ci\delta.
$
Let us, however, provide also a direct proof.
Since $A_c$ and $B_c^*$ generate  $C(S^2_{qc})$ in the $C^*$-algebraic sense, and both 
$\bar{r}_1^c$ and $\delta$ are continuous, it suffices to check the desired
equality on $A_c$ and $B_c^*$.
Taking advantage of (\ref{ra}) and using the fact that $\delta$ is a 
$C^*$-homomorphism,
we obtain
\[
(\delta\ci \bar{r}_1^c)(A_c)
=\delta(f_c(A))=f_c(\delta(A_c))=f_c(A_c)\ot 1 =\bar{r}_1^c(A_c)
\ot 1=
((\bar{r}_1^c\ot id)\ci\delta)(A_c).
\]
To handle $B_c^*$ 
it is useful to consider
its polar decomposition $B_c^*=V_c|B_c^*|$ (see (\ref{pol})).
Now, $\delta(B_c^*)=B_c^*\ot{u^*}^2$ entails 
$\delta(B_cB_c^*)=B_cB_c^*\ot u^2{u^*}^2=|B_c^*|^2\ot 1$, 
whence, by the continuity of the square root function, 
$\delta(|B_c^*|)=|B_c^*|\ot 1$. Consequently,
\[
\delta(V_c)(|B_c^*|\ot 1)=\delta(B_c^*)=
B_c^*\ot{u^*}^2=(V\otimes {u^*}^2)(|B_c^*|\ot 1).
\]
Thus, due to the invertibility of $|B_c^*|\ot 1$, we have $\delta(V_c)=
V_c\ot {u^*}^2$.
 (We chose $B_c^*$ rather than $B_c$ because, unlike $|B_c|$, $|B_c^*|$ is
invertible.) 
On the other hand, $\bar{r}_1^c(V_c)=V_c$ (see above) and $\bar{r}_1^c(|B_c^*|)=
g_c(A_c)$
(see (\ref{rb*})). To complete the proof, one can reason in the same way as for
generator $A_c$.
\ere

We end this section by recalling $K$-facts for the quantum disc.
Since $C(D^2_{q^4})$ is isomorphic to the Toeplitz algebra,
the ``standard" exact sequence~\cite[p.68]{b-b98} 
(cf.\ (\ref{cd}), (\ref{qs}), (\ref{qrp})) is equivalent to: 
\beq \Label{kdi}
0\lra{\cal K}\lra C(D^2_{q^4})
   \lra C(S^1)\lra 0, 
\eeq
from which it follows that 
\beq\Label{kdi2}
K_0 (C(D^2_{q^4}))\cong\IZ,~~~K_1(C(D^2_{q^4}))\cong 0~~~
\mbox{\cite[p.123]{w-ne93}}.
\eeq

\section{Quantum real projective space }
\setcounter{equation}{0}

Our first aim is to define on the equator quantum sphere $S^2_{q\infty}$ a
$\Z_2$-action mimicking the antipodal action of $\Z_2$ on $S^2$. 
The geometrical meaning of generators (see \ref{cart2}) hints at the
formulas $r_2(A)=-A$, $r_2(B)=-B$. Owing to the even nature of algebraic
relations in $P(S^2_{q\infty})$, these equalities indeed define the desired
action on $P(S^2_{q\infty})$. (Note that this recipe would not work for 
$P(S^2_{qc})$, $c\in[0,\infty)$.)
The $*$-algebra of quantum real projective 2-space
can now be defined by 
\bde[\cite{h-pm96}]
$P(\IR P^2_q)=\{a\in P(S^2_{q\infty})|r_2(a)=a\}$. 
\ede
\bre\em
Recall that $\dr\ci r_2=(r_2\ot id)\ci\dr$
(see above Section~6 in \cite{p-p87}), where, much as before,
 $\dr$ is the  restriction to $P(S^2_{q\infty})$
of the coproduct $\hD$ in $P(SU_q(2))$.
(Since both $r_2$ and $\dr$ are algebra homomorphisms, it suffices
to check this formula on generators, where it is evidently true.)
Thus the antipodal action and the $SU_q(2)$-action on the equator quantum
sphere are compatible.
Consequently, just as quantum spheres themselves,
 $\R P^2_q$ is an (embeddable) quantum homogeneous space of $SU_q(2)$,
i.e., $\hD(P(\R P^2_q))\inc P(\R P^2_q)\ot P(SU_q(2))$.
\ere
Unlike the quantum disc, $\IR P^2_q$ is a $\Z _2$-Galois quotient
of the equator quantum sphere, i.e., $S^2_{q\infty}\ra\IR P^2_q$
is an (algebraic) quantum principal bundle \cite[Proposition 2.10]{h-pm96}.
 As mentioned in the proof of Proposition~\ref{main}, 
the elements $A^kB^l$, $A^m{B^*}^n$, $n>0$, form a basis
of $P(S^2_{q\infty})$. Taking this into account,
 it is straightforward that $P(\IR P^2_q)$ is the 
$*$-subalgebra of $P(S^2_{q\infty})$ generated by $A^2,~B^2$ and $AB$.
We put 
\beq\Label{prt}
P=A^2,~~R=B^2,~~T=AB,
\eeq
and find immediately the following relations:
\beq\Label{ppp}
P=P^*,~~~ T^2=q^2PR,~~~ RT^*=q^2T(-q^4P+I),~~~ R^*T=q^{-2}T^*(-P+I),
\eeq
\beq\Label{pcpr}
RR^*=q^{12}P^2-q^4(1+q^4)P+I,~~~ R^*R=q^{-4}P^2-(1+q^{-4})P+I,
\eeq
\beq\Label{pcpt}
TT^*=-q^4P^2+P,~~~ T^*T=q^{-4}(P-P^2).
\eeq
\beq\Label{pc}
RP=q^8PR,~~~ RT=q^4TR,~~~ PT=q^{-4}TP,
\eeq
\bpr\Label{uni}
Let $I_q$ be the $*$-ideal in the free $*$-algebra $\IC\<P,R,T\>$ generated by 
the relations (\ref{ppp})-(\ref{pc}).
Then the $*$-algebra $\IC\<P,R,T\>/I_q$ is isomorphic to $P(\IR P^2_q)$.
\epr
\bpf
There exists a $*$-algebra epimorphism $f:\IC\<P,R,T\>/I_q\ra P(\IR P^2_q)$
given on generators by $f(P)=A^2,~f(R)=B^2,~f(T)=AB.$ On the other hand,
we can define a linear map $g:P(\IR P^2_q)\ra \IC\<P,R,T\>/I_q$ by its 
values on the elements of a basis of $P(\IR P^2_q)$:
$g(A^{2k}B^{2l})=P^kR^l,~~g(A^{2k+1}B^{2l+1})=P^kTR^l,~~
g(A^{2m}{B^*}^{2(n+1)})=P^m{R^*}^{n+1},~~ 
g(A^{2m+1}{B^*}^{2n+1})=q^2P^mT^*{R^*}^n$.
Evidently, $f\ci g=id$. Consequently $g$ is injective and the above elements of 
$\IC\<P,R,T\>/I_q$ are linearly independent. To have the reverse equality 
$g\ci f=id$ it suffices to show that these elements span $\IC\<P,R,T\>/I_q$.
Assume inductively that every monomial in $P,R,T,R^*,T^*$
 of length at most $n$ is in the span. 
This is clearly true for $n=1$. Take now an arbitrary monomial $M_{n+1}$ of
length $n+1$. It can always be written as $M_nW$, where $M_n$ is a 
monomial of length $n$ and $W$ is one of the elements $P,~R,~T,~R^*,~T^*$.
By assumption $M_n$ is a linear combination of 
$P^kR^l,~P^kTR^l,~P^m{R^*}^{l+1},~P^mT^*{R^*}^l$. Using the commutation relations
(\ref{ppp})-(\ref{pc}) among generators, it can be directly verified 
that each of the monomials $M_nW$ is again in the span.
\epf

In order to extend the antipodal $\Z_2$-action to $C(S^2_{q\infty})$, 
note first that (\ref{delta2}) entails
\beq
{r}_2(A)=(\bar{r}_1\ci\delta_{\sqrt{-1}})(A),~~~
{r}_2(B)=(\bar{r}_1\ci\delta_{\sqrt{-1}})(B).
\eeq
Therefore, we can define the completion of $r_2$ by 
\beq\Label{rtwo}
\bar{r}_2:=\bar{r}_1\ci\delta_{\sqrt{-1}} : C(S^2_{q\infty})\lra C(S^2_{q\infty}).
\eeq
Observe that we need to put $g=\sqrt{-1}$
rather than $g=-1$  because this 
$U(1)$-action comes from $SU(2)$ which is the double-cover of $SO(3)$.
Therefore, to rotate the quantum sphere by the angle $\pi$ 
(antipodal action is such a rotation composed with reflection),
we take $g=e^{i\pi/2}$ rather than $g=e^{i\pi}$.
Since both $\bar{r}_1$ and $\delta_{\sqrt{-1}}$ are
$C^*$-homomorphisms, we can define the $C^*$-algebra of $\R P^2_q$ as
\bde\Label{trp2}
$
C(\R P^2_q):=\{a\in C(S^2_{q\infty})~|~
\bar{r}_2(a)=a\}.
$
\ede
 Arguing as in Proposition~\ref{co} (second
part of the proof), we get
that the completion of $P(\R P^2_q)$ with respect to the norm on 
$C(S^2_{q\infty})$ coincides with the thus defined $C(\R P^2_q)$. 
To study the structure of this $C^*$-algebra, let us prove:
\bth\Label{t2}
There are no unbounded $*$-representations of the $*$-algebra
$P(\R P^2_q)$. Up to the unitary equivalence, all irreducible
(bounded) $*$-representations of this algebra  are the following:

(i) A family of one-dimensional representations 
$\rho_\theta:P(\R P^2_q)\ra \C$ parameterized by $\theta\in[0,2\pi)$,
which are given by
\beq
\rho_\theta(P)=\rho_\theta(T)=0,~~\rho_\theta(R)=e^{i\theta}.
\eeq

(ii) An infinite dimensional representation $\rho$ (in a Hilbert
space $H$ with an orthonormal basis $\{e_k\}_{k\in\N}$) given by

\beq
\rho(P)e_k=q^{4k}e_k,
\eeq
\beq
\rho(T)e_k=\left\{\ba{cc}0&k=0\\q^{2(k-1)}(1-q^{4k})^{1/2}e_{k-1}&k>1
\ea\right .,
\eeq
\beq
\rho(T^*)e_k=q^{2k}(1-q^{4(k+1)})^{1/2}e_{k+1},~~k\geq 0,
\eeq
\beq
\rho(R)e_k=\left\{\ba{cc}0&k=0,1\\(1-q^{4k})^{1/2}(1-q^{4(k-1)})^{1/2}e_{k-2}&
k>1\ea\right .,
\eeq
\beq
\rho(R^*)e_k=(1-q^{4(k+1)})^{1/2}(1-q^{4(k+2)})^{1/2}e_{k+2},~~k\geq 0.
\eeq
\ethe
\bpf
Suppose that $\rho$ is an unbounded $*$-representation 
(\cite[Definition 8.1.9]{s-k90}) of $P(\R P^2_q)$. 
The relations $T^*T=q^{-4}(P-P^2)$ and $P=P^*$ entail that
both $P$ and $1-P$ are positive. Thus we have $0\leq\rho(P)\leq 1$,
so that  $\rho(P)$ is bounded. It follows then from the relations (\ref{pcpr})
 and (\ref{pcpt}) that also $\rho(R)$ and
$\rho(T)$ are bounded.
Therefore, unlike $P(D^2_{q^4})$, the $*$-algebra $P(\R P^2_q)$
  has no unbounded representations.

Now, let $\rho$ be an irreducible bounded 
$*$-representation in a Hilbert space $H$.
As before, we have $0\leq\rho(P)\leq 1$.
Let $\rho(P)=0$. Then $T^*T=q^{-4}(P-P^2)$ implies that also 
$\rho(T)=0$. Hence the only remaining relation is $\rho(R^*R)=1=\rho(RR^*)$,
and we can see that the image of $\rho$ is commutative.
Since the only irreducible representations of
a commutative algebra are one-dimensional, we arrive at~(i).

Let us now assume $\rho(P)\neq 0$. It is immediate from $RP=q^8PR$ and
$PT=q^{-4}TP$ that $\ker\rho(P)$ is $\rho$-invariant. Due to the irreducibility
and the boundedness of $\rho$, either $\ker\rho(P)=H$ or $\ker\rho(P)=0$. Since 
the first case is
excluded by the assumption $\rho(P)\neq 0$, we have $\ker\rho(P)=0$.
Using the characterization of elements of the spectrum by
approximate eigenvectors 
and taking advantage of the relation $PT=q^{-4}TP$ we will
show that the spectrum of $\rho(P)$ consists of the eigenvalues
$q^{4k},~k\in\N$, and their limiting point $0$.
We already know that the spectrum of $\rho(P)$ lies in the interval
$[0,1]$. Next, note that $0$ cannot be the only element of $Sp(\rho(P))$ because
this would mean $\rho(P)=0$, contradicting $\ker(\rho(P))=0$. For the same reason,
$0$ cannot be an eigenvalue. If $1$ would be the only element of the
spectrum, we would have $\rho(P)=1$, and consequently, due to 
$PT=q^{-4}TP$, $\rho(T)$ would vanish. This would contradict $\rho(TT^*)=1-q^4$ 
resulting from the relation $TT^*=-q^4P^2+P$. Thus $1$ cannot be the only
element in the spectrum. It is also impossible that $Sp(\rho(P))=\{0,1\}$
because then 0 would be an eigenvalue.
Summing up, we have shown that there exists 
$\lambda\in Sp(\rho(P))\cap (0,1)$, and that $\rho(T)\neq 0$.

By \cite[Lemma 3.2.13, vol.1]{kr97}, there exists a sequence $(\xi_n)_{n\in \N}$
of unit vectors in the representation space $H$ such that
\beq\Label{app}
\lim_{n\rightarrow\infty}\|\rho(P)\xi_n-\hl\xi_n\|=0.
\eeq
We will now show that there exist $N\in\N$ and $C>0$ such that 
$\|\rho(T)\xi_n\|\geq C$
for $n\geq N$. To estimate $\|\rho(T)\xi_n\|$, we use $T^*T=q^{-4}(P-P^2)$.
Now, the right hand side of this equality we want to put in a form allowing us
to apply (\ref{app}). Adding and subtracting
$\hl^2-\hl$ gives:
\beq
(P-P^2)= (P-\hl) + (\hl^2-P^2) - (\hl^2-\hl)
= (1-\hl-P)(P-\hl)+ \hl(1-\hl).
\eeq
Therefore, using the triangle inequality and $\|a\|\|\eta\|\geq\|a\eta\|$,
we obtain
\beq\Label{r1}
\|(\rho(P)-\rho(P^2))\xi_n\|
\geq |\hl(1-\hl)| - \|1-\hl-\rho(P)\| \|(\rho(P)-\hl))\xi_n\|.
\eeq
On the other hand,
\[\Label{r2}
\|\rho(T^*)\|\|\rho(T)\xi_n\|\geq\|\rho(T^*T)\xi_n\|=
q^{-4}\|(\rho(P)-\rho(P^2))\xi_n\|.
\]
Combining (\ref{r1}) with (\ref{r2}) and remembering that 
$\|\rho(T)^*\|=\|\rho(T)\|\neq 0$, we get
\beq\Label{absch}
\|\rho(T)\xi_n\|\geq\frac{|\hl(1-\hl)|}{q^{4}\|\rho(T^*)\|}
-\|(\rho(P)-\hl)\xi_n\|
\frac{\|1-\hl-\rho(P)\|}{q^4\|\rho(T^*)\|}.
\eeq
Since $\frac{|\lambda(1-\lambda)|}{q^{4}\|\rho(T^*)\|}$ is positive,
the existence of the desired $N$ and $C$ follows from~(\ref{app}). 
Hence we conclude that
\[
\eta_n:=\frac{\rho(T)\xi_n}{\|\rho(T)\xi_n\|}
\]
are well-defined unit vectors for $n\geq N$.
Our goal now is to show
\beq\Label{q-4l}
\lim_{n\rightarrow\infty}\|\rho(P)\eta_n-q^{-4}\lambda\eta_n\|=0,
\eeq
which is tantamount to $q^{-4}\hl\in Sp(\rho(P))$. Assume $n\geq N$.
Then $\|\rho(T)\xi_n\|\geq C$.
Using $PT=q^{-4}TP$, we have
\beq
\|\rho(P)\eta_n-q^{-4}\lambda\eta_n\|
=\frac{\|\rho(T)\rho(P)\xi_n-\hl\rho(T)\xi_n\|}{q^{4}\|\rho(T)\xi_n\|}
\leq \frac{\|\rho(T)\|}{q^{4}C}
\|\rho(P)\xi_n-\hl\xi_n\|.
\eeq
Consequently, (\ref{q-4l}) follows from (\ref{app}). Thus we have shown that
\[\Label{ci}
\hl\in Sp(\rho(P))\cap(0,1)\Rightarrow q^{-4}\hl\in Sp(\rho(P)).
\]
 Therefore,
there exists $k$ such that $q^{-4k}\hl=1$, for otherwise we would get an
unbounded sequence $\hl,q^{-4}\hl,q^{-8}\hl,\ldots,q^{-4k}\hl,\ldots\in Sp(\rho(P))$
contradicting $Sp(\rho(P))\inc[0,1]$. Hence 
$Sp(\rho(P))$ $\subseteq\{q^{4k}~|~k\in\N\}\cup\{0\}$. It also follows that
$1\in Sp(\rho(P))$.
Thus 1 is isolated in $Sp(\rho(P))$, so that it is an eigenvalue, and there exists
a vector $\xi$ such that $\rho(P)\xi=\xi$, $\|\xi\|=1$.
(Notice that now we evidently have $\|\rho(P)\|=1$.)

It follows from the relation $T^*P=q^{-4}PT^*$ that 
$\rho({T^*}^k)\xi$ are eigenvectors of $\rho(P)$ corresponding 
to the eigenvalue $q^{4k}$.
Let us prove inductively that all these eigenvectors are different from zero.
For $n=0$, the statement $\rho({T^*}^n)\xi \neq 0$ is automatically true.
Assume now that $\rho({T^*}^n)\xi \neq 0$ for some $n\in\N$. Then,
using $TT^*=-q^4P^2+P$ and $T^*P=q^{-4}PT^*$, one obtains
\[\Label{tn1}
T{T^*}^{n+1}={T^*}^n(q^{4n}P-q^{8n+4}P^2).
\]
 Therefore,
$
\rho(T)\rho({T^*}^{n+1})\xi=(q^{4n}-q^{8n+4})\rho({T^*}^n)\xi\neq 0
$
 and
consequently $\rho({T^*}^{n+1})\xi\neq 0$.
Hence, by induction, $\rho({T^*}^{n})\xi\neq 0,~\fa n\in\N$.
This proves that  $Sp(\rho(P))=
\{q^{4k}~|~k\in\N\}\cup\{0\}$. Since $\rho(P)$ is self-adjoint and 
$\rho({T^*}^k)\xi$ are eigenvectors 
of different
eigenvalues of $\rho(P)$, they are mutually orthogonal. Thus, the vectors
\beq\Label{ek}
e_k:=\frac{\rho({T^*}^k)\xi}{\|\rho({T^*}^k)\xi\|},~~~ k\in\N,
\eeq
form an orthonormal system. 
Let us now show that the span of the $e_k$'s is closed under the 
action of the entire algebra.
We already know that 
\[
\rho(P)e_k=q^{4k}e_k.
\]
On the other hand, the formula (\ref{tn1}) entails
\bea\Label{txi}
\|\rho({T^*}^{k+1})\xi\|
&=&
\langle\rho(T)\rho({T^*}^{k+1})\xi,\rho({T^*}^k)\xi\rangle^{1/2}
\nonumber\\
&=&
\langle\rho({T^*}^k)(q^{4k}-q^{8k+4})\xi,\rho({T^*}^k)\xi\rangle^{1/2}
\nonumber\\
&=&
q^{2k}(1-q^{4(k+1)})^{1/2}\|\rho({T^*}^k)\xi\|.
\eea
Hence, from the definition (\ref{ek}) we have
\[
\rho(T^*)e_k=q^{2k}(1-q^{4(k+1)})^{1/2}e_{k+1}.
\]
 The relation 
$T^*T=q^{-4}(P-P^2)$
implies that $\rho(T)e_0$ has zero length: 
\[
\|\rho(T)e_0\|^2=\<e_0,\rho(T^*)\rho(T)e_0\>=0.
\]
 Thus $\rho(T)e_0=0$.
Similarly, 
$R^*R=q^{-4}P^2-(1+q^{-4})P+1$ entails that
$\rho(R)e_0=\rho(R)e_1=0$.
Using the equality $TT^*=q^4P^2+P$ one obtains
\[
\rho(T)e_k=q^{2(k-1)}(1-q^{4k})^{1/2}e_{k-1},~ k>0.
\] 
 Furthermore, a straightforward computation
taking advantage of
$RT^*=q^2T(-q^4P+1)$ and $RT=q^4TR$ gives
\[
\rho(R)e_k=(1-q^{4k})^{1/2}(1-q^{4(k-1)})^{1/2}e_{k-2},~ k\geq 2,
\]
and
\[
\rho(R^*)e_k=(1-q^{4(k+1)})^{1/2}(1-q^{4(k+2)})^{1/2}e_{k+2},~~k\geq 0.
\]
Therefore the Hilbert space $H_e:=\overline{span\{e_k\}}$ 
is a closed invariant subspace of $H$, 
and we have $H_e=H$ by the irreducibility of the bounded representation $\rho$.
Any other irreducible representation $\rho'$ with $\rho'(P)\neq 0$  generates 
an orthonormal basis in the same way, so it has to be
unitarily equivalent to the above one.
\epf
Observe that the irreducible representations of the above theorem 
are restrictions of representations of $C(S^2_{q\infty})$. 
More precisely, we have that $\rho_\theta$ is the restriction of $\pi_\theta$, and
 $\rho$ is the restriction of $\pi_+$. Since every representation in a
 separable Hilbert space is a direct
integral of irreducible representations, we can conclude that all representations
of $C(\R P^2_q)$ extend to representations of $C(S^2_{q\infty})$. Therefore,
the norm of the universal $C^*$-algebra of $P(\R P^2_q)$ coincides with the norm
inherited from $C(S^2_{q\infty})$. Note also \footnote{
We owe this observation to P.~Podle\'s.}  that $\rho$ is a restriction of 
both $\pi_+$ and $U\pi_-U^{-1}$, where $U$ is the unitary defined by 
$Ue_k=(-1)^ke_k$.
Now the faithfulness of $\rho$ follows 
from the faithfulness of $\pi_+\oplus\pi_-$.
Again, since the norm in any representation is
always less or equal to the norm in a faithful representation,
we can conclude that
the universal and inherited norms coincide.
Summarizing we have established:
\bco
The $C^*$-algebra $C(\R P^2_q)$ is the {\em universal} $C^*$-algebra of $P(\R P^2_q)$.
\eco

\bre\em
Similarly to the case of the reflection action $\bar{r}^c_1$, we want the diagram
\vspace*{.3cm}
\begin{center}
\xext=2000 \yext=500
\begin{picture}(\xext,\yext)
\setsqparms[1`1`1`1;1000`500]
\putsquare(0,0)[C(S^2_{qc})`C^*(\Sh)\oplus_\sigma C^*(\Sh)`
C(S^2_{qc})`C^*(\Sh)\oplus_\sigma C^*(\Sh);
\pi_+^c\oplus\pi_-^c`\bar{r}_2^c` `\pi_+^c\oplus\pi_-^c]
\setsqparms[-1`1`1`-1;1000`500]
\putsquare(1000,0)[\phantom{C^*(\Sh)\oplus_\sigma C^*(\Sh)}`C(S^2_{q\infty})
`\phantom{C^*(\Sh)\oplus_\sigma C^*(\Sh)}`C(S^2_{q\infty});
\pi_+\oplus\pi_-` `\bar{r}_2`\pi_+\oplus\pi_-]
\end{picture}
\end{center}
\vspace*{.3cm}
to be commutative. To this end, we define the antipodal  
$\Z_2$-actions on the equilateral spheres by
$
\bar{r}_2^c:=\chi_c^{-1}\ci\bar{r}_2\ci\chi_c.
$ 
It is clear that
 $\{a\in C(S^2_{qc})|\bar{r}_2^c
(a)=a\}\cong C(\R P^2_q)$. 
Furthermore, it follows directly from  definitions, Proposition~\ref{com}
and (\ref{coma}) that 
$\bar{r}_2^c
=\bar{r}_1^c\ci\hd_{\sqrt{-1}}
=\hd_{\sqrt{-1}}\ci\bar{r}_1^c$.
Remembering also that $U(1)$ is Abelian, we have that the antipodal
actions are compatible with the $U(1)$-actions on quantum spheres:
$
\hd_g\ci\bar{r}_2^c=\bar{r}_2^c\ci\hd_g.
$
\ere

Let us turn now to the computation of $K$-groups of $C({\R}P^2_q)$.
Just as in the classical case (e.g., 
see \cite[Corollary 6.47]{k-m78}), we have: 
\bth\Label{krp}
The topological K-groups of the quantum real projective space
$\R P^2_q$ are as follows:
$K_0(C({\R}P^2_q))\cong {\Z}\oplus{\Z}_2$,  
$K_1(C({\R}P^2_q)) \cong 0$.
\ethe
\bpf
First we need to find an exact sequence analogous to~(\ref{crp}). 
Let $J$ be the closed two-sided $*$-ideal of $C({\R}P^2_q)$ generated 
by $P$, and let $p:C({\R}P^2_q)\rightarrow C({\R}P^2_q)/J$ 
be the natural surjection. 
Arguing as in the proof of Theorem~\ref{t2}, we see that
in the quotient all relations 
(see (\ref{ppp})-(\ref{pc})) reduce to $p(RR^*)=I=p(R^*R)$.
 Consequently $p(R)$ 
is unitary, and we have $C({\R}P^2_q)/J\cong 
C^*(p(R))\cong C(S^1)$. 
\ble\Label{idj}
The ideal $J$ is isomorphic (via the faithful representation
$\rho$) to the $C^*$-algebra $\cal K$ of 
compact operators on a separable Hilbert space.
\ele
\bpf
The operator $\rho(P)$ is evidently compact 
(see Theorem~\ref{t2} for an explicit formula), 
whence $\rho(J)\inc {\cal K}$. On the other
hand, as $\rho(P)$ is a diagonal operator with 
eigenvalues of multiplicity one,  all the one-dimensional projections $P_k$
onto the vectors $e_k$ are elements of $\rho(J)$.  Furthermore, since also
$\rho(T)$ belongs to $\rho(J)$ (see above the lemma) and it
 is a weighted shift with non-vanishing coefficients,
all matrix units $E_{ij}$ belong to $\rho(J)$. 
(They can be obtained from the $P_k$ and 
$\rho(T)$.) Therefore, ${\cal K}\inc\rho(J)$, and consequently ${\cal K}
=\rho(J)$.
 The claim of the lemma
follows from the faithfulness of $\rho$.
\epf
Denote by $i$ the inverse of the appropriate restriction of $\rho$,
and again by $p$ the canonical surjection $p$ composed with the isomorphism
$C({\R}P^2_q)/J\cong C(S^1)$. With the help of Lemma \ref{idj}, we obtain
the desired exact sequence:  
\beq
\Label{qrp}
 0\lra{\cal K}\stackrel{i}{\lra}
   C({\R}P^2_q)\stackrel{p}{\lra}C(S^1)\lra 0. 
\eeq
Since $K_1({\cal K})\cong 0$, $K_1(C(S^1))\cong\Z\cong K_0({\cal K})$,
the 6-term exact sequence of $K$-theory yields 
\[\Label{6} 
0\lra K_1(C({\R}P^2_q))\lra{\Z}
   \stackrel{\partial}{\lra}{\Z}
   \stackrel{i_*}{\lra}K_0(C({\R}P^2_q))
   \stackrel{p_*}{\lra}{\Z}\lra 0, 
\] 
where $\partial:{\Z}\cong K_1(C(S^1))\rightarrow 
K_0({\cal K})\cong{\Z}$ is the index map. 
Due to the exactness of the above sequence, to compute the $K$-groups
it suffices to determine the index map~$\partial$.
We know from the preceding discussion that $p(R)$ is the unitary
generator of $C(S^1)$. Hence $[p(R)]$ generates $K_1(C(S^1))$. 
We have $[p(R)]\cong 1$ via the identification of $K_1(C(S^1))$ with~$\Z$.
Thus all we need to complete the calculation is the value of $\partial$ on $[p(R)]$.
In general, if $A$ is a unital $C^*$-algebra, 
 $0\ra I\ra A\ra A/I\ra 0$ the
short exact sequence inducing the 6-term exact sequence, 
$u$  a unitary element of $A/I$, and $\nu\in A$ 
such that $\nu\nu^*=1$ and $\nu/I=u$, then $\partial([u])=[1-\nu^*\nu]$
(see \cite[Section 8.3.2]{b-b98} or \cite[Remark 8.1.4]{w-ne93}).
We need to lift the unitary $p(R)$ to an appropriate coisometry
in $C({\R}P^2_q)$. 
As $\rho(R)$ is a weighted double-shift, 
the desired coisometry could be given by $U(e_n)=e_{n-2}$, $Ue_{1}=0=Ue_{2}$. 
The operator $U$
satisfies the polar decomposition $\rho(R)=U|\rho(R)|$.
Furthermore, since $p(|R|)=1$ and $\rho$ is faithful, 
if there exists $\nu_2\in C({\R}P^2_q)$ such that $\rho(\nu_2)=U$,
then $p(R)=p(\nu_2)$ and $\nu_2$ is the desired coisometry. 
Thus we need to show that $U\in\rho(C({\R}P^2_q))$.
Let ${\cal U}$ be an open neighbourhood of
$\{1,q^4\}$ such that ${\cal U}\cap \{q^{4k}\}_{k=2,...,\infty}=\emptyset$,
and let $\hb:{\R}\rightarrow{\R}$ be a continuous function given on 
$\R\setminus{\cal U}$ by the formula
\[
\hb(x)=
(1-x)^{-1/2}(1-q^{-4}x)^{-1/2}.
\]
 Then $U=\rho(R)\hb(\rho(P))$, whence
$U\in\rho(C({\R}P^2_q))$, as needed. Consequently, 
\[
\partial([p(R)])=[1-\nu_2^*\nu_2]=[diag(1,1,0,...)]\cong 2.
\]
Here in the penultimate equality we identified $K_0(C({\R}P^2_q))$ with
$K_0(\rho(C({\R}P^2_q)))$, and in the last step
$K_0({\cal K})$ with $\Z$. 
Finally, as $\partial$ is injective, (\ref{6}) breaks into two exact sequences:
\[
0\lra  K_1(C({\R}P^2_q))\lra 0
\]
and
\[\Label{es3}
0\lra{\Z}
   \stackrel{2\cdot}{\lra}{\Z}
   \stackrel{i_*}{\lra}K_0(C({\R}P^2_q))
   \stackrel{p_*}{\lra}{\Z}\lra 0.
\]
The former gives immediately $K_1(C({\R}P^2_q))\cong 0$, and the latter splits,
as $\Z$ is a free module over itself. Therefore we conclude that 
$K_0(C({\R}P^2_q))\cong \mbox{Im}~i_*\oplus\Z$. On the other hand, the exactness
of (\ref{es3}) implies the exactness of
\[
0\lra 2\Z\lra\Z\st{i_*}{\lra}\mbox{Im}~i_*\lra 0.
\]
Hence $\mbox{Im}~i_*\cong \Z_2$, and consequently 
$K_0(C({\R}P^2_q))\cong\Z_2\oplus\Z$.
\epf
\vspace*{-1cm}
\bre\em
From the exact sequence (\ref{6})
we can read the generators of $K_0(C({\R}
P^2_q))$: $[I]$ of infinite order and $[f]$ of order 
2, where $f$ is a minimal projection in the ideal $J$.    
\ere
It is an immediate consequence of Theorem~\ref{krp} and (\ref{kdi2})
that $C({\R}P^2_q)$ is {\em not} the standard extension of $C(S^1)$
by $\cal K$:
\bco
The $C^*$-algebras $C({\R}P^2_q)$ and $C(D^2_{q^4})$ are {\em not}
isomorphic, i.e.,  $C({\R}P^2_q)$ is not the Toeplitz algebra.
\eco
\bre\Label{Jones}\em
The $C^*$-algebra $C(S^2_{q\infty})$ is generated by $A$ and $B$. 
The antipodal ${\Bbb Z}_2$-action sends $A$ to $-A$, $B$ to $-B$, 
and $C({\Bbb R}P^2_q)$ is the fixed-point subalgebra. 
There is a conditional expectation 
(e.g., see \cite[pp.570,571 in vol.2]{kr97}) 
$E:C(S^2_{q\infty})\rightarrow C({\Bbb R}P^2_q)$ 
sending even monomials to themselves and annihilating the odd ones.
We have a ``quasi-basis" $\{u_1:=I,~u_2:=A,~u_3:=B^*\}$ for $E$ 
(in the sense of \cite[p.2]{w-y90}).
 This means that 
\[
a=\sum_{i=1}^3 E(au_i)u_i^*=\sum_{i=1}^3 u_iE(u_i^*a),~
\fa a\in C(S^2_{q\infty}).
\]
Hence 
\[
{\rm Index}(E):=\sum_{i=1}^3 u_iu_i^*=I+A^2+B^*B=2I. 
\]
We can think of $S^2_{q\infty}$ as a two-fold covering of~${\Bbb R}P^2_q$.
\ere

{\bf Acknowledgments.}
The authors are indebted to M.~Bo\.{z}ejko, D.~Calow, L.~D\c{a}browski, P.~Podle\'s, 
W.~Pusz, K.~Schm\"{u}dgen, J.C.~Varilly, S.L.~Woronowicz and J.~Wysocza\'nski
for very helpful discussions.
This work
was partially supported by the KBN grant 2 P03A 030 14 and the 
Naturwissenschaftlich-Theoretisches Zentrum of Leipzig University. 
P.M.H.\ and R.M.\ are grateful for hospitality to Leipzig University,  
Max Planck Institute for Mathematics in the Sciences, and Warsaw
University, Polish Academy of Sciences, respectively.


\end{document}